\documentclass[11pt]{amsart}

\usepackage{amssymb}
\usepackage{amsmath}
\usepackage{enumerate}
\usepackage{amstext}
\usepackage{graphicx}
\usepackage{hhline}
\usepackage{psfrag}
\usepackage{float}
\usepackage[dvips]{epsfig}
\usepackage{rotating}

\newtheorem{theorem}{Theorem}[section]

\newtheorem{remark}[theorem]{Remark}
\newtheorem{rem}[theorem]{Remark}

\DeclareSymbolFont{msbm}{U}{msb}{m}{n}
\DeclareMathSymbol{\C}{\mathalpha}{msbm}{'103}
\DeclareMathSymbol{\R}{\mathalpha}{msbm}{'122}
\DeclareMathSymbol{\Z}{\mathalpha}{msbm}{'132}
\DeclareMathSymbol{\N}{\mathalpha}{msbm}{'116}

\title{A class of multi-phase traffic theories for  microscopic, kinetic and continuum traffic models}

\author{Raul Borsche \and Mark Kimathi
        \and Axel Klar  }

\address{
Raul Borsche,  Mathematics Department,
University of Kaiserslautern \newline P.O.Box 3049, 67653
Kaiserslautern, Germany.
\newline {\rm \texttt{Email: borsche@mathematik.uni-kl.de}
}
\newline
Mark Kimathi,  Mathematics Department,
University of Kaiserslautern \newline P.O.Box 3049, 67653
Kaiserslautern, Germany.
\newline {\rm \texttt{Email: kimathi@mathematik.uni-kl.de}
}
\newline
Axel Klar, Mathematics Department,
University of Kaiserslautern, \newline P.O.Box 3049, 67653
Kaiserslautern, Germany; Fraunhofer ITWM Kaiserslautern, 67663 Kaiserslautern, Germany
\newline {\rm \texttt{Email: klar@mathematik.uni-kl.de}
}
}

\begin{document}

\keywords{ traffic flow, macroscopic equations, kinetic
derivation, multi-valued fundamental diagram, stop and go waves, phase
transitions}

\begin{abstract}
In the present paper  a review and numerical comparison of  a special class of multi-phase traffic theories based on  microscopic, kinetic and macroscopic traffic models is given.
Macroscopic  traffic equations with
multi-valued fundamental diagrams are derived  from different microscopic and kinetic models.
Numerical experiments show similarities and differences of the models, in particular, for  the appearance and structure of stop and go
waves for highway traffic in dense situations.
For all   models, but one, phase transitions can appear near bottlenecks depending on the
local density and velocity of the flow.  
\end{abstract}

\maketitle
\markboth{Raul Borsche,Mark Kimathi
        and Axel Klar}{Continuum multi-phase traffic theories}

\section{Introduction}
\label{Introduction}

Traffic flow modeling has been considered on different levels of description, see \cite{BD11} for a recent review:
on the microscopic level the motion of each   vehicle is described.
Mathematical models are generally stated using a  large system of ordinary differential equations for  position and velocity  of the  vehicles
based on Newtonian mechanics
 \cite{BD99,May90, BHNSS95,GHR61,HB01}.
On the macroscopic level  the state of the system  is described by
averaged quantities. Typically, density and  linear momentum are used to describe the flow. 
The corresponding mathematical models are based on systems of nonlinear partial differential equations derived from conservation laws with suitable closure relations.
Starting from the pioneering  work of  Aw and
Rascle \cite{AR98}  new macroscopic models for traffic flow 
have been derived and investigated intensively in the last decade, see for example \cite{Deg, CG07,AKMR03,Gre00,Goa06,Hel01}. 
These models avoid several inconsistencies of previous models, like wrong way traffic and missing bounds on the density.
We note that these models can be derived from microscopic models  in a variety of ways, see for example \cite{AKMR03,Zha}.
Finally, kinetic theory describes  the state of the system 
by a probability distribution function of the   position and velocity of the vehicles. Mathematical
models generally use
integro-differential or Fokker-Planck type equations.
Kinetic equations for vehicular traffic can be found, for example, in
\cite{PH71,PF75,Nel95,KW97}. Procedures to derive macroscopic traffic
equations including the Aw/Rascle model from underlying kinetic models
have been performed in different ways by several authors, see, for
example, \cite{Hel95B} and \cite{KW00}. These procedures are developed
in analogy to the transition from the kinetic theory of gases to
continuum gas dynamics.

Another  basic problem of macroscopic traffic
flow equations has been described by Kerner \cite{Ker98,Ker99,Ker00}.
The observations there suggest a more complicated dependence of the
homogeneous steady speed states on density: these states are not given
by a uniquely defined function $u= U^e (\rho)$ as in the above mentioned 
models, but cover a whole range in the density-flow diagram leading to a multi-valued fundamental diagram.
The resulting dynamical system has a multi-phase behavior  in the sense of Kerner: the flow changes between different  stationary state which represent free and so called synchronized or jam behavior.
In the context of the derivation of macroscopic models from microscopic ones  the homogeneous steady state solutions can be interpreted as  an emergent behavior of interactions at the microscopic scale and  multiple solutions may be related to the heterogeneous behavior of the driver-vehicle subsystem.
A variety of microscopic and macroscopic multiphase models has been developed by several authors.
In particular, there is a large number of works on microsopic models. We refer among many others to
\cite{Ker03,Ker04,RKP11,SH09,TKH10,ZK05,BM09,THH00}.
Macroscopic  models for traffic flow with   phase transitions in the sense of Kerner  can be found in  \cite{BWGPB,Goa06,Col2, Col}. 
However, these models 
do not describe phenomena like stop and go waves near bottlenecks.
For microscopic and macroscopic multiphase models exhibiting    stop and go waves and similar traffic instabilities we refer  to  \cite{Ker03, Ker04, ZK05, LL12, GKR03,SM06}.
Kinetic models  allowing for multiple stationary solutions and associated macroscopic models with multi-valued fundamental diagrams
 have been developed in \cite{IH,IKM,GKMW03, NS}.
We refer to \cite{BD11}   for  a recent review and to \cite{DCB99,BDC02, Dag97}  for  further material on the above issues.

The present paper contains a comparison and discussion of a class of  macroscopic models with multi-valued fundamental diagrams.
We consider  models of the form
 \begin{eqnarray}
\partial_t \rho + \partial_x (\rho u ) &=&0 ,\\
 \partial_t (\rho u)  + \partial_x ( \rho u^2 )  - c(\rho)  \partial_x u &=& \rho R ( u,   \tau) \nonumber
\end{eqnarray}
 with right hand side 
 \begin{eqnarray*}
R (\rho,u)  = \; \frac{\rho}{T} \;
\left[U(\rho,u) -  u \right] 
\end{eqnarray*}
and fundamental diagrams given by functions $U= U(\rho,u)$ having at least two  equilibrium solutions, i.e. solutions of the equation
$u = U(\rho,u)$ for fixed $\rho$ out of a  certain density domain, where multiphase traffic may appear.

The paper starts with a review of this class of macroscopic multi-phase models.
The models are derived from either microscopic or kinetic equations. To guarantee a proper comparison of the considered models 
several changes to the original models are proposed. 
Moreover, the parameters of the different models are chosen  such that the stable equilibrium solutions are the same for all models.
Then, the  different models are numerically  investigated for a bottleneck problem 
and the appearance of stable wave patterns is shown which can be interpreted as stop-and go waves at (on ramp) bottlenecks.
This numerical comparison as well as the changes made for each of the models  to make them comparable are new up to the knowledge of the authors.
We note that  stable periodic waves excited by small periodic perturbations have been studied in a series of papers also for 
equations with single valued right hand sides, see \cite{Gre04, GKR03,FKNRS09}. Remarks on these waves can be found in section
\ref{Numerical Investigations}.

The paper is arranged in the following way: In Section \ref{A continuum  multi-phase traffic model}
 the derivation of macroscopic equations from microscopic models is reviewed and applied to a multi-phase traffic model  from \cite{Ker00,Ker03}.
 In Section \ref{Multi-phase hydrodynamic equations derived from kinetic equations}  kinetic equations are investigated
 and used to derive  multi-valued fundamental
diagrams.  The different models are partially changed to make them comparable to each other. Section \ref{Comparison of multi-phase hydrodynamic models} contains a summary and comparison of the different approaches and the derived multi-valued fundamental diagrams.
  Finally, in Section
\ref{Numerical Investigations} numerical results are given comparing  the
different density-velocity relations.  Moreover, an inhomogeneous traffic
flow situation with a bottleneck is investigated, showing the
appearance of traffic instabilities together with a qualitative comparison of the structure of these instabilities.

\section{Continuum  multi-phase traffic model derived from microscopic equations}
\label{A continuum  multi-phase traffic model}

\subsection{From microscopic to macroscopic models}
\label{micromacro}
We review the classical procedure for so called 'General Motors' (GM) type car-following models, see \cite{BD99,May90}.
 Denoting with 
 $x_i(t), v_i(t), i = 1,
\ldots, N$  the  location and
speed of the vehicles  at time $t \in \R^+$, 
and the distance between successive cars by
 \begin{eqnarray*}
l_i &=& x_{i+1} - x_{i},
\end{eqnarray*}
we consider the microscopic  equations 
\begin{eqnarray*}
\dot{x}_i &=& v_i\\
\dot{v}_i &=&  
C \frac{(v_{i+1} -v_i)}{l_{i}-H}
+  \frac{1}{T} (U(\rho_i) -v_i).
\nonumber
\end{eqnarray*}
The local 
``density around vehicle i'' and its inverse 
(the local (normalized) ``specific volume'') are respectively 
defined by
\begin{eqnarray*}
\rho_i = \frac{H}{l_i} \; \mbox{and} \;\; 
\tau_i = \frac{1}{\rho_i} = \frac{l_i}{H},
\end{eqnarray*}
where $H$ is the length of a car.

\begin{rem}
\label{rem1}
Here,  the density is normalized and 
therefore dimensionless, so that 
the {\em maximal} density 
is 
$\rho_m = 1/\tau_m = 1$.
\end{rem}
The constant $C>0$ 
and the relaxation time $T$ 
are given parameters. The function $U= U(\rho), 0 \le  \rho \le \rho_m =1$ is the so called fundamental diagram, see \cite{PH71, Ker04, SM06}.
The simplest choice is given by $U(\rho) = 1-\rho$.
One obtains
the  microscopic model
\begin{eqnarray}
\label{ode2}
\dot{x}_i &=& v_i \;,\\
\dot{v}_i &=& \frac{C}{H} \; \frac{(v_{i+1} -v_i) }
{\tau_{i}-1}
+   \; \frac{1}{T} \left(U(\rho_i) -v_i\right) . \nonumber
\end{eqnarray}
We have
\begin{eqnarray*}
\dot{l}_i = v_{i+1}-v_i
\;\; \mbox{or} \;\; 
\dot{\tau}_i = \frac{1}{H} \; (v_{i+1}-v_i).
\end{eqnarray*}

The limit of number of cars going to infinity yields the Lagrangian form of the macroscopic equations, see \cite{AKMR03}.
We obtain the equivalent of the p-system  in gas dynamics (isentropic Euler equations in Lagrangian form), compare \cite{GR},
\begin{eqnarray}
\label{pde}
\partial_T \tau - \partial_X u &=& 0 \; ,\\
\partial_T u - c (\rho)  \partial_X u &=&  \;\frac{1}{T} \left[U(\rho) -  u \right] \; ,
\nonumber
\end{eqnarray}
where $\tau = \frac{1}{\rho}$ is the specific volume, i.e. the (local) 
dimensionless fraction of space occupied by the cars.
 $\rho$  the (normalized) density is the limit of $\rho_i$ defined above, as the number of cars tends to infinity.
 $u$ is  the 
macroscopic 
velocity of the cars. 
Moreover, 
\begin{equation}\label{prho}
c (\rho)=
C  \left( \frac{1}{\rho} -1 \right)^{-1}
\end{equation}
and the  function
$
U(\rho) 
$
is defined in the microscopic model above. 
We change the Lagrangian
``mass'' coordinates $(X,T)$ into Eulerian coordinates  $(x,t)$ with
$$\partial_x X= \rho , \;\partial_t X = - \rho v , \;
T = t$$
or
$$\partial_X x= \rho^{-1} = \tau,  \;
\partial_T x =  v.$$
Thus, $X= \int^x \rho(y,t) dy$ 
 describes 
the 
total space occupied by cars up to point $x$. 
The macroscopic 
system in Eulerian coordinates is then 
\begin{eqnarray}
\label{cons}
\partial_t \rho + \partial_x (\rho u ) &=&0 ,\\
 \partial_t (\rho u)  + \partial_x ( \rho u^2 ) - c(\rho) \partial_x u  &=& \; \frac{\rho}{T} \;
\left[U(\rho) -  u \right] .\nonumber
\end{eqnarray}

For the well-posedness of the above problem
we refer to \cite{Ber}.

\begin{remark}
The same approach works for right hand sides with multi-valued equilibrium distributions
\begin{eqnarray}
\; \frac{\rho}{T} \;
\left[U(\rho,u) -  u \right] .
\end{eqnarray}

Examples are Switching curve (SC) models as in  \cite{GKR03} with

\begin{eqnarray*}
U ( \rho,u) =  
\begin{cases}
U_1(\rho) ,
&   \rho < \rho_{f}  \; \mbox{or}  \;  u > S(\rho),  \rho_{f}< \rho<\rho_{j} ,\\
U_2(\rho) ,
&   u < S(\rho),  \rho_{f}< \rho<\rho_{j}  \mbox{or}  \;\rho > \rho_{j}.
\end{cases}
\end{eqnarray*}
Here the switching curve $S= S(\rho)$ is given. We consider  the density  in  the (synchronized flow)  region between a lower bound of free flow $\rho_f$  and an upper bound of jam traffic $\rho_j$. Then, there exists multiple stationary states $U_1, U_2$ whose regions of influence are separated by the switch-curve $S$.
 
 A similar model is the  Speed-adaptation (SA)  type models of Kerner et.al  \cite{Ker03} with
 \begin{eqnarray*}
 U (   \rho,u) =  
 \begin{cases}
 U_1(\rho ) ,  & u > U_{sync}   , \rho < \rho_{j} ,\\
 U_2(\rho),
 &  u < U_{sync}  ,\rho >  \rho_{f} ,\\
\end{cases}
 \end{eqnarray*}
where the parameter $U_{syn}$ is the averaged speed,  which separates the domains of influence of the two stationary states
in the 2-D region of synchronized flow  in the flow-density plane.
 \end{remark}

\subsection{A microscopic ATD-type  model}

In this section we sketch a  simplified version of the microscopic Acceleration time delay (ATD) model of Kerner et. al \cite{Ker03}.
We consider a microscopic model with the variables space, velocity and acceleration:

\begin{eqnarray}
\label{multiphasemicro}
\dot{x}_i &=& v_i\\
\dot{v}_i &=&  a_i \nonumber\\
\dot{a}_i &=& F(a_i, v_i, \frac{v_{i+1} - v_i}{H} , \frac{x_{i+1} - x_i}{H} )
\nonumber
\end{eqnarray}

with
\begin{eqnarray*}
F(a_i, v_i, \frac{\Delta v_{i} }{H} , \tau_i ) =  
\begin{cases}
(a_i^{free}- a_i)/T_{del} ,
& \tau_i > G (v_i)   , \tau_i > \tau_{jam}  ,\\
(a_i^{syn}- a_i)/T_{del} ,
&  \tau_i < G (v_i) ,\tau_i > \tau_{jam} ,\\
(a_i^{jam}- a_i)/T_{del} ,
&  \tau_i < \tau_{jam} .
 \\\end{cases}
\end{eqnarray*}
Here,  $\Delta v_{i}  = v_{i+1} - v_i$ and $a_i^{free}, a_i^{syn}, a_i^{jam}$   denote the desired accelerations in the free, synchronized and jam region respectively and $T_{del}$ denotes the time delay of the acceleration of the vehicle.   The function $G$ separates the free from the synchronized acceleration behaviour  and will  be fixed later at the end of Section \ref{Comparison of multi-phase hydrodynamic models}.
For  a proper comparison with the above models we change the definitions in\cite{Ker03} of the different accelerations slightly and define
the terms as  follows
\begin{eqnarray*}
a_i^{free} &=& \frac{1}{T} (U(\rho_i ) - v_i) + \frac{ c(\rho_i)}{H} \Delta v_i \\
a_i^{syn}&=& \frac{1}{T} \min (U(\rho_i ) - v_i,0 ) +  \frac{ c(\rho_i)}{H}  \Delta v_i\\
a_i^{jam}&=&  -   \frac{1}{T} v_i+ \frac{ c(\rho_i)}{H}  \Delta v_i,
\end{eqnarray*}
where
$
U(\rho)$
and $c(\rho)$ are given  as before. This means acceleration depends  on the speed difference to the predecessor and a term relaxing to a desired acceleration.

\subsubsection*{The hydrodynamic multi-phase model}

 To obtain the  hydrodynamic version of the microscopic model in the last section we follow the procedure in Section \ref{micromacro}.
 In Lagrangian coordinates we obtain directly
 \begin{eqnarray}
\label{pdemulit}
\partial_T \tau - \partial_X u &=& 0 \; ,\\
\partial_T u &= & a \nonumber \\
\partial_T a &=&  F(a,u,\partial_X u , \tau)  .
\nonumber
\end{eqnarray} 
This leads to the following equations in Eulerian coordinates
 \begin{eqnarray}
\label{consmulti}
\partial_t \rho + \partial_x (\rho u) &=&0 ,\\
 \partial_t (\rho u)  + \partial_x ( \rho u^2 )  &=& \rho a \nonumber\\
\partial_t (\rho a)  + \partial_x ( \rho u  a)  &=& \rho F (a,u,\tau \partial_x u , \tau) . \nonumber
\end{eqnarray}
 
\subsubsection*{A reduced model}

 Assuming that the delay times for acceleration are small we can reduce the above ATD-type model to 
  \begin{eqnarray}
\label{consred}
\partial_t \rho + \partial_x (\rho u ) &=&0 ,\\
 \partial_t (\rho u)  + \partial_x ( \rho u^2 )  - c(\rho)  \partial_x u &=& \rho R ( u,   \tau), \nonumber
\end{eqnarray}
 where 
 \begin{eqnarray*}
R( \rho,u) =  
\begin{cases}
\frac{1}{T}  (U(\rho ) - u),  & \tau > G (u)   , \tau > \tau_{j} ,\\
\frac{1}{T}  \min (U(\rho) - u,0 )  ,
&  \tau < G (u)  ,\tau> \tau_{j} ,\\
-  \frac{1}{T} u,
&  \tau < \tau_{j} 
 .\\\end{cases}.
\end{eqnarray*}

This is equivalent to 
\begin{eqnarray*}
R (\rho,u)  = \; \frac{\rho}{T} \;
\left[U(\rho,u) -  u \right] 
\end{eqnarray*}

with
 \begin{eqnarray*}
U( \rho, u) =  
\begin{cases}
U(\rho ),  & \tau > G(u)  , \rho < \rho_{j} \;  \mbox{or } \; \tau < G(u) , u> U(\rho), \rho < \rho_{j} ,\\
u ,
&   \tau < G(u)  , u < U(\rho), \rho< \rho_{j} ,\\
0,
&  \rho > \rho_{j} 
 .\\\end{cases}
\end{eqnarray*}
 
For comparison with the other multi-valued fundamental diagrams we rewrite the relaxation term  using $K(u) = 1/G(u)$:
\begin{eqnarray*}
U ( \rho,u) =  
\begin{cases}
U(\rho ) ,  & \; \rho <  K(u)  , \rho < \rho_{j} \;  \mbox{or } \; \rho > K(u) , u> U(\rho), \rho_f <  \rho < \rho_{j} ,\\
u ,
&  \rho > K(u)  , u < U(\rho), \rho_{f}< \rho<\rho_{j} ,\\
0 ,
&  \rho > \rho_{j} 
 .\\\end{cases}
\end{eqnarray*}
 
\section{Multi-phase hydrodynamic equations derived from kinetic equations}
\label{Multi-phase hydrodynamic equations derived from kinetic equations}

\subsection{Kinetic equations and correlations}

The basic quantity in a kinetic approach is the single car
distribution $f (x,v)$ describing the density of cars at $x$ with
velocity $v$.  The total density $\rho$ on the highway is
\begin{eqnarray*}
 \rho (x) \;=\; \int_0^w  f(x,v) dv,
\end{eqnarray*}
where $w$ denotes the maximal velocity.  Let $F(x,v)$ denote the
probability distribution in $v$ of cars at $x$, i.e. $f(x,v) = \rho(x)
F(x,v)$.  The mean velocity is
\begin{eqnarray*}
 u (x) \;=\; \int_0^w v  F(x,v) dv.
\end{eqnarray*}

An important role is played by the distribution $f^{(2)} (x,v,h,v_+)$
of pairs of cars being at the spatial point $x$ with velocity $v$ and
leading cars at $x+h$ with velocity $v_+$.  This distribution function
has to be approximated by the one-vehicle distribution function
$f (x,v)$. We use the chaos assumption
\begin{eqnarray*}
f^{(2)} (x,v,h,v_+) = q(h)\, f(x,v)\, F(x+h,v_+),
\end{eqnarray*}
compare  \cite{Nel95}.  For a vehicle with velocity $v$ the
function $q(h,v;\rho,u)$ denotes the distribution of leading vehicles
with distance $h$ under the assumption that the velocities of the
vehicles are distributed according to the distribution function $f$.

Thresholds for braking ($H_{B}$)
and acceleration ($H_{A}$) are introduced. 
From a microscopic point of view drivers will brake,
once the distance between the driver and its leading car is becoming
smaller than a threshold $H_B$ and will accelerate, once this distance
is becoming larger than $H_A$. Otherwise the cars will not change the
velocities.  Velocities are changed instantaneously once acceleration
or braking lines are reached.  Models including acceleration of the
cars can be developed as well, see \cite{IKM} for an example.

The distribution of leading vehicles $q(h)$ is prescribed a
priori.  The main properties, which $q(h)$ has to fulfill are
positivity, 
\begin{eqnarray*}
\int_0^\infty q(h)dh \;=\; 1, 
\end{eqnarray*}
and
\begin{eqnarray}
\label{meanvalue}
 \int_0^{\infty} h q(h) dh  \;=\;  \frac{1}{\rho}.
\end{eqnarray}
Equation (\ref{meanvalue})
 means that the average headway of the cars is $1/\rho$.
The leading vehicles are assumed to be distributed in an uncorrelated
way with a minimal distance $H_B$ from the car under
consideration, see  \cite{Nel95}:
\begin{eqnarray*}
q(h) \;=\;
\tilde{\rho}\, e^{-\tilde{\rho}(h-H_B)}\,\chi_{[H_B,\infty)}(h).
\end{eqnarray*}
The reduced density $\tilde\rho$ has to be defined in such a way, that
(\ref{meanvalue}) is fulfilled. One obtains
\begin{eqnarray}
\label{admissible}
\tilde\rho\;=
\;=\;\frac{\rho}{1-\rho\, H_B}.
\end{eqnarray}

We note that 
\begin{eqnarray*}
q_A = q(H_A) 
\;=\; \tilde\rho\,e^{-\tilde\rho(H_A-H_B)}
\end{eqnarray*}
and
\begin{eqnarray*}
q_B = q(H_B) 
\;=\; \tilde\rho\;.
\end{eqnarray*}

The probability $P_{ov} =P_{ov} (\rho,u)$ for overtaking or lane changing
and the
corresponding probability $P_B = 1- P_{ov}$ for braking
are determined from 
phenomenological  considerations:
at constant density, free flow of cars, i.e. larger velocities  will be related to larger probabilites of overtaking or smaller probabilites of braking.
So called synchronized traffic is associated to smaller velocities 
and thus larger  probabilites of braking.
That means the probablitiy of braking can be considered as  - for fixed density -- a decaying function of velocity $u$.
Similar  arguments can be found for example in \cite{Ker00}.

\begin{remark}
In the following we present a kinetic model. Note that the results
like multi-valued fundamental diagrams and stop and go behaviour of
the derived macroscopic equations do not depend on the
exact choice of the microscopic interactions we have  chosen here.
The model discussed in the
next section is only chosen due to the fact that explicit
stationary solutions are available. We could as well have chosen models like in \cite{KW00} or Fokker-Planck type models like in \cite{IKM}.
\end{remark}

\subsection{The evolution equation}

To write the kinetic evolution equations in a simple form we use
\begin{eqnarray*}
k \;=\; k(\rho,u) \;=\; \frac{P_B\,q_B}{q_A + P_B q_B} 
\end{eqnarray*}
and
\begin{eqnarray*}
\gamma \;=\; \gamma (\rho,u) \;=\; \frac{q_A }{1-k} \;=\; q_A + P_B q_B.
\end{eqnarray*}
We consider a relaxation frequency 
\begin{eqnarray*}
\nu \;=\; \nu (k)
\end{eqnarray*}
and define 
\begin{eqnarray*}
\frac{1}{T} = \gamma\,\nu .
\end{eqnarray*}

The kinetic  model is then given by the following evolution equation
for the distribution function $f$:

\begin{eqnarray}
\label{kinetic2}
\partial_t  f + v \partial_x  f
&=&  C^+ (f)\\
&=&
\gamma \Big[ k (G^+_B -L^+_B ) (f) + (1-k)(G^+_{A} -L^+_{A}) (f) \nonumber  \\
&&+ \nu  (G_{S} -L_{S}) (f) \Big]  \nonumber
\end{eqnarray}    
with the loss and gain terms for braking interactions
\begin{eqnarray*}
G_{B}^+ (f)
&=&
\int \int_{\hat{v} > \hat{v}_+}
\vert \hat{v} -\hat{v}_+ \vert
\sigma_{B} (v ; \hat{v},\hat{v}_+) 
 f(x,\hat{v})
 F(x+H_B, \hat{v}_+)
d \hat{v} d\hat{v}_+\\
L_{B}^+ (f)
&=&
\int_{\hat{v}_+<v} 
\vert v -\hat{v}_+ \vert
f(x,v)
 F(x+H_B, \hat{v}_+) d\hat{v}_+.
 \end{eqnarray*}  
 The loss and gain terms for acceleration is defined as 
 \begin{eqnarray*}
 G_{A}^+ (f)&=&
\int \int_{\hat{v} < \hat{v}_+}
 \vert \hat{v} -\hat{v}_+ \vert
\sigma_{A} (v ; \hat{v}, \hat{v}_+) 
 f(x,\hat{v})
 F(x+H_A, \hat{v}_+)
d \hat{v} d\hat{v}_+\\
L_{A}^+ (f)&=&
\int_{\hat{v}_+>v}  
\vert v -\hat{v}_+ \vert 
 f(x,v)
 F(x+H_A, \hat{v}_+)  d\hat{v}_+.
  \end{eqnarray*}  
 Finally  terms  describing the random behavior of drivers are 
   \begin{eqnarray*}  G_{S} (f)&=&
\int_{0}^w
\sigma_{S} (v , \hat{v})
f(x, \hat{v})
d \hat{v}\\
\vphantom{\int_{0}^w}
L_{S} (f)&= & f(v).
\end{eqnarray*}

 $\sigma_B$ and $\sigma_A$ denote the  distribution  of the new velocities $v$ after an 
interaction.   Reaching the braking line the vehicle
brakes, such that the new velocity $v$ is distributed with a
distribution function $\sigma_B$ depending on the old velocities
$\hat{v},\hat{v}_+$.
For acceleration the new velocity is  distributed according to $\sigma_A$.
 The  relaxation term is
introduced to include  a random behaviour of the drivers. 

\begin{remark}
For further details  on this Boltzmann/Enskog approach to traffic flow
modelling see \cite{KW97}.  
\end{remark}

\subsection*{Example}

For the  probability distributions $\sigma_A, \sigma_B$
we  choose the following simple expressions:
\begin{eqnarray}
\label{1s}
\sigma_{B} (v , \hat{v},\hat{v}_+) \;=\; \frac{1}{\hat{v} - \hat{v}_+}
\,\chi_{[\hat{v}_+, \hat{v}]} (v)
\end{eqnarray}
and
\begin{eqnarray}
\label{2s}
\sigma_{A} (v , \hat{v} , \hat{v}_+) \;=\; \frac{1}{\hat{v}_+ -  \hat{v}}
\,\chi_{[\hat{v},  \hat{v}_+]} (v).
\end{eqnarray}
This means we have an equidistribution of the new velocities
between the velocity of the car and the velocity of its leading
car. Finally,
\begin{eqnarray}
\label{3s}
\sigma_{S} (v , \hat{v}) \;=\; \frac{1}{w} .
\end{eqnarray}

\subsection{Stationary Distributions 
and multi-valued Fundamental Diagrams}
\label{Stationary Distributions 
and multi-valued Fundamental Diagrams}

In this section we investigate the stationary homogeneous equations
and determine the multi-valued fundamental diagrams.
We consider the local interaction operator:
\begin{eqnarray*}
C(f) \;=\; \gamma \left[k (G_B -L_B ) (f) + (1-k)(G_{A} -L_{A}) (f)
+ \nu (G_{S} -L_{S}) (f) \right] 
\end{eqnarray*}
with $f = \rho F$. The gain and loss terms $G_B, L_B$, etc. are
defined in the same way as $G_B^+, L_B^+$, etc., except that  $x+ H_X ,
X=A,B$ is substituted by $x$, wherever it appears.  The homogeneous
stationary equation is
\begin{eqnarray*}
C(f) \;=\;  0.
\end{eqnarray*}       
We assume that for fixed $\rho$ and  $k$ there is a
unique solution
\begin{eqnarray*}
f \;=\; f^e \;=\; \rho F^e (k,v)
\end{eqnarray*}
of this equation. This is true for the example stated above.

Thus, for  fixed $k$ the mean value of $F^e$ is then
\begin{eqnarray*}
u^e(k) \;:=\; \int_0^w v F^e(k,v) dv.
\end{eqnarray*}
The function
$u^e$ is uniquely determined due to the above assumption as a 
function of $k$.  However, this does not yield immediately the 
fundamental diagram, i.e. an equilibrium relation between flux and density.

Instead, the
fundamental diagram  is determined  from the following considerations:
let  $u$  be  the (possibly multi-valued) solution of the equation
\begin{eqnarray}
\label{multiu}
  u  &\;=\;  u^e(k(\rho,u))
\end{eqnarray}
for fixed $\rho$. 
If there is a unique solution we obtain a well defined
relation for equilibrium velocity and density and the usual fundamental 
diagram.
However, in general this equation will have a
multitude of different solutions $u$, even infinitely many.  Plotting
a dependence of this solution on the density one obtains in the
general case a two-dimensional region in the density-velocity plane, where
the solutions are located.
The fundamental diagram is then a multi-valued 
function.

\begin{remark}

For the example above the homogeneous solution can be solved explicitly  and the corresponding multi-valued solutions of equation (\ref{multiu}) can be evaluated numerically. Explicit expressions for $F^e(k)$  and $u^e(k) $ can be found in \cite{GKMW03}.
A plot of $u^e(k)$ is shown in 
Figure \ref{uek}.

\end{remark}

\begin{remark}
In contrast to the other models described above the kinetic approach gives an explanation for a multi-valued  fundamental diagram using in particular the braking probability $P_B$ as a basic quantity.
\end{remark}

\subsection{Derivation of  Macroscopic  Models}
\label{Derivation of Macroscopic Models}

In this section macroscopic
equations for density and mean velocity are derived.
Different procedures are described, for example in \cite{KW00,IP09}.

\subsection*{Balance Equations}
\label{Balance Equations}

Multiplying the inhomogeneous kinetic equation (\ref{kinetic2}) with
$1$ and $v$ and integrating it with respect to $v$ one obtains the
following set of balance equations:
\begin{eqnarray}
\label{balance}
\partial_t \rho + \partial_x (\rho u) & = & 0\\
\partial_t (\rho u) + \partial_x ( P + \rho u^2)
+ E
& = & S\nonumber
\end{eqnarray}
with the 'traffic pressure'
\begin{eqnarray}
\label{pressure}
P \;=\;\int_0^w  (v-u)^2 f dv,
\end{eqnarray}
the  Enskog flux term
\begin{eqnarray}
\label{ens}
E \;=\; \int_0^w v [C(f)(x,v,t) - C^+(f) (x,v,t)] dv,
\end{eqnarray}
and the source term
\begin{eqnarray}
\label{source}
S \;=\; \int_0^w  v C( f)(x,v,t) dv.
\end{eqnarray}

For the present discussion we are, in particular, interested in  the source term S.

\subsection*{Closure }
\label{Closure Relations}

We concentrate on the relaxation term and cite the results for the other terms , compare \cite{KW00}.
The traffic
pressure $P$ is negligible and approximated by zero, see \cite{KW00}.
Moreover, the Enskog term $E$ is approximated by linearizing
expression (\ref{ens}) for $E$ in $H$. 
We obtain \cite{KW00}
\begin{eqnarray*}
E \;\sim\; -   c_{kin}  (\rho) \partial_x u
\end{eqnarray*}
with $c_{kin} (\rho)$ given by the details of the collision operator. In the following we will neglect the special form of $c_{kin}$
and choose $c_{kin}= c(\rho)$ for comparison as in the other models described above.

Finally, the source term $S$ has to be approximated. We use a relaxation time approximation
\begin{eqnarray*}
C (f)  \;\sim\;   \frac{1}{T} \left(f^e(k (\rho,u),v) -f(v)\right).
\end{eqnarray*}
This yields 
\begin{eqnarray*}
S \;\sim\; S^{e}(\rho,u) \;=\; \rho  \frac{1}{T} \left(u^e(k (\rho,u)) -u\right).
\end{eqnarray*}

Thus, from the kinetic approach one obtains macroscopic equations of the form
\begin{eqnarray}
\label{macro}
\partial_t \rho + \partial_x (\rho u) &=& 0\\
\partial_t (\rho u ) + \partial_x(\rho u^2)
-   c (\rho)  \partial_x u
&=& S^{e} (\rho,u) \nonumber
\end{eqnarray}
with
\begin{eqnarray*}
 S^{e}(\rho,u) \;=\; \rho \frac{1}{T} \left(u^e(k (\rho,u)) -u\right),
\end{eqnarray*}
where  $k = k(\rho,u)$ is defined as 
\begin{eqnarray*}
k \;=\; \frac{P_B\,q_B}{q_A + P_B q_B}  =  \frac{1}{1+ \frac{\mbox{exp}(- \tilde \rho (H_A - H_B)}{P_B}} . \end{eqnarray*}
Choosing $H_A= H_B$ this simplifies to 
\begin{eqnarray*}
k \;=\;  \frac{1}{1+ \frac{1}{P_B}} . \end{eqnarray*}

\begin{remark}
One obtains a multi-valued variant of the Aw-Rascle equations with
a multi-valued relaxation term on the right hand side.
\end{remark}

\section{Comparison of multi-phase hydrodynamic models}
\label{Comparison of multi-phase hydrodynamic models}

We consider  models of the form
 \begin{eqnarray}
\partial_t \rho + \partial_x (\rho u ) &=&0 ,\\
 \partial_t (\rho u)  + \partial_x ( \rho u^2 )  - c(\rho)  \partial_x u &=& \rho R ( u,   \tau) \nonumber
\end{eqnarray}
 with
 \begin{equation}
c (\rho)=
C  \left( \frac{1}{\rho} -1 \right)^{-1}
\end{equation}
and 
 \begin{eqnarray*}
R (\rho,u)  = \; \frac{\rho}{T} \;
\left[U(\rho,u) -  u \right] 
\end{eqnarray*}
and fundamental diagrams given by functions $U$ of the following form:

\subsection*{ATD-type models}

\begin{eqnarray*}
U ( \rho,u) =  
\begin{cases}
U(\rho ) ,  & \; \rho <  K(u)  , \rho < \rho_{j} \;  \mbox{or } \; \rho > K(u) , u> U(\rho), \rho_f <  \rho < \rho_{j} ,\\
u ,
&  \rho > K(u)  , u < U(\rho), \rho_{f}< \rho<\rho_{j} ,\\
0 ,
&  \rho > \rho_{j} 
 .\\\end{cases}
\end{eqnarray*}

\subsection*{SA-type models}

 \begin{eqnarray*}
 U (   \rho,u) =  
 \begin{cases}
 U_1(\rho ) ,  & u > U_{sync}   , \rho < \rho_{j} ,\\
 U_2(\rho)
 &  u < U_{sync}  ,\rho >  \rho_{f} ,\\
\end{cases}
 \end{eqnarray*}

See Kerner et.al. \cite{Ker03} for details.

\subsection*{Switching curve models}

\begin{eqnarray*}
U ( \rho,u) =  
\begin{cases}
U_1(\rho) ,
&   \rho < \rho_{f}  \; \mbox{or}  \;  u > S(\rho),  \rho_{f}< \rho<\rho_{j} ,\\
U_2(\rho) ,
&   u < S(\rho),  \rho_{f}< \rho<\rho_{j}  \; \mbox{or}  \;\rho > \rho_{j}.
\end{cases}
\end{eqnarray*}
where $S(\rho)$  is a switching curve. For an investigation of these models, see \cite{GKR03}.

\subsection*{Kinetic models:}

$$U (\rho,u) =  u^e (k(\rho,u))   $$
with 
\begin{eqnarray*}
k \;=\;  \frac{1}{1+ \frac{1}{P_B}} .\end{eqnarray*}

In Figure \ref{uekinetic} we plot the  equilibrium  solutions of $ u =  U (\rho,u)$ together with  the values $ U (\rho,u) -u $ denoting the length of the arrows.  Moreover,  $U (\bar \rho,u) -u$ is plotted for fixed $\rho$ versus $u$ in Figures \ref{uekineticcut}.

For a proper comparison of the above models  the parameters are chosen as follows:

\begin{eqnarray*}
U(\rho) = U_1(\rho)\\
K^{-1} (\rho)  = U_2(\rho),   \rho_j > \rho >\rho_f  \\
K(u) = \rho_f,    u> U_2 (\rho_f)\\
U_2 (\rho)<U_{sync}  < U_1 (\rho),  \rho_f < \rho < \rho_j.
\end{eqnarray*}

Moreover, the functions $\nu$ and $P_B$ in the kinetic model are chosen such that the stable kinetic equilibrium solutions of $u = u^e(k(\rho,u))$ are given by $U_1$ and $U_2$, the instable solution by $S(\rho)$.
This leads to the fundamental diagrams shown in Figures \ref{uekinetic} to \ref{ueATD}.

\section{Numerical Investigations}
\label{Numerical Investigations}
 
In this section we compare the different  approaches numerically.
First, the kinetic model is investigated and the associated fundamental diagram is determined.
Second the other three models stated in the last section are compared to the kinetic approach
and third all three approaches are used in a inhomogeneous traffic simulation with a bottleneck.

\subsection{The  stationary, homogeneous kinetic equation}

We consider the kinetic equation and resulting fundamental diagrams.
For the numerical simulations we normalize and use  $w=1$.

\begin{figure}[ht]
\begin{center}
\epsfig{file=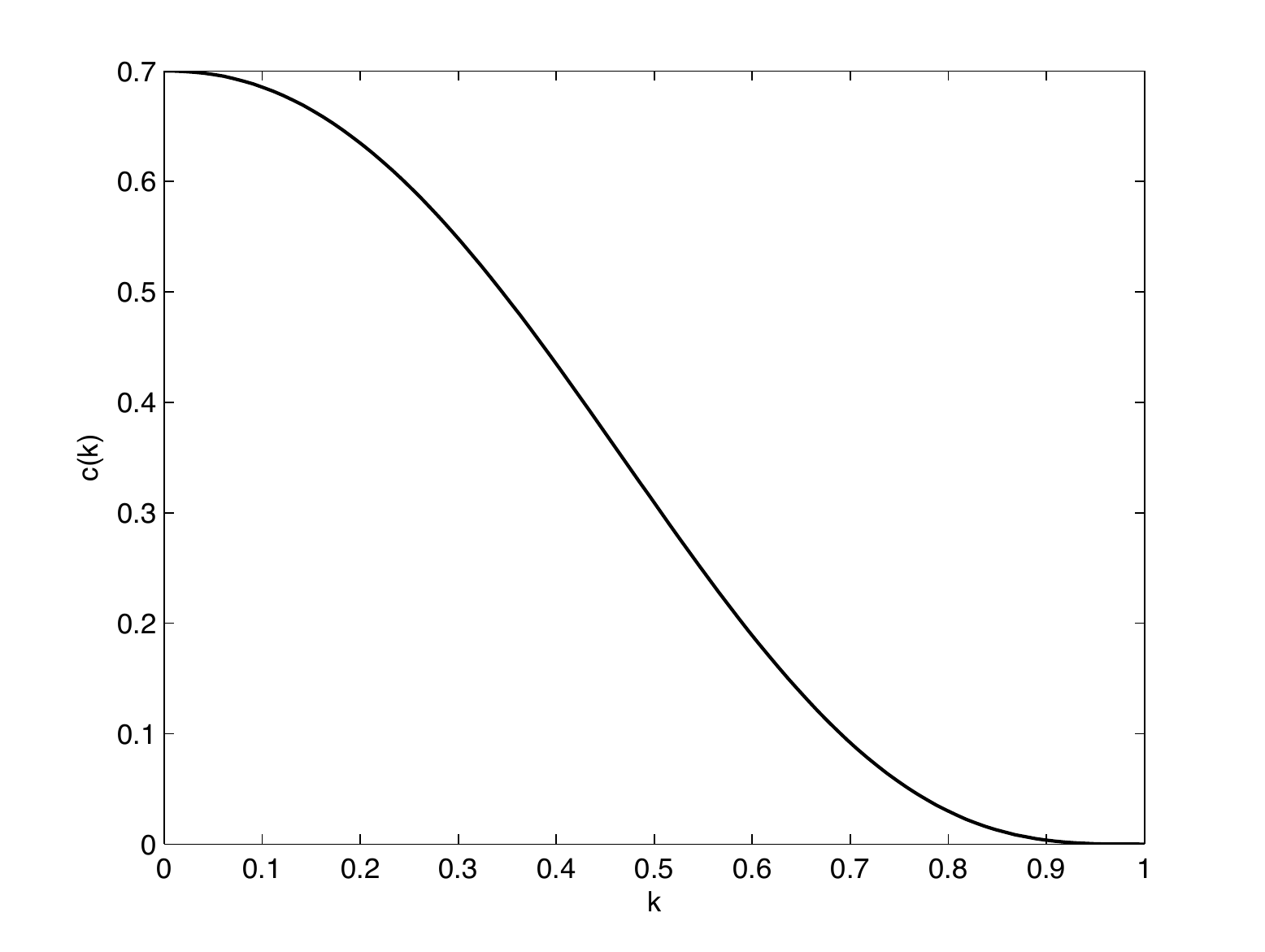,height=.25\textheight,width=.45\textwidth}
\epsfig{file=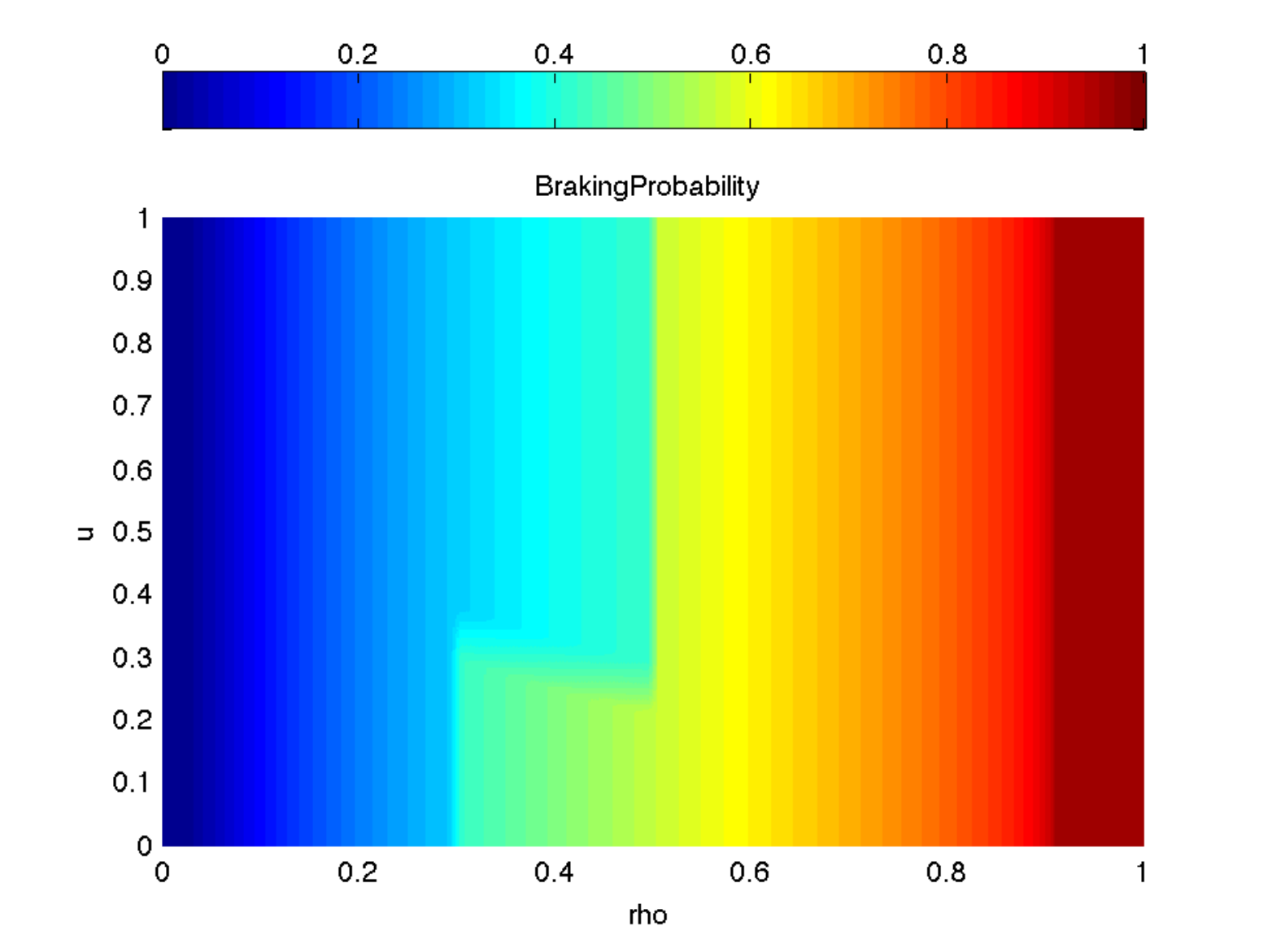,height=.25\textheight,width=.45\textwidth}
\end{center}
\caption{ Left: Frequency of random events $\nu(k)$. Right: Braking Probability $P_B (\rho,u)$.}
\label{cNpb}
\end{figure}

Moreover, we  choose
$\nu$ as in Figure \ref{cNpb}.
A reasonable function
$\nu$ should be zero for maximal density ($k=1$).
In this  case  there is no more random behaviour of the drivers,
all drivers have velocity $0$.
For the case $k=0$ we have chosen $\nu$ as a finite quantity. 
If  these two features are fulfilled, the qualitative behaviour of the model
does not depend on the exact form of $\nu$.
The  braking probability $P_B$ is plotted as well in Figure \ref{cNpb}.

\begin{figure}
\begin{center}
\epsfig{file=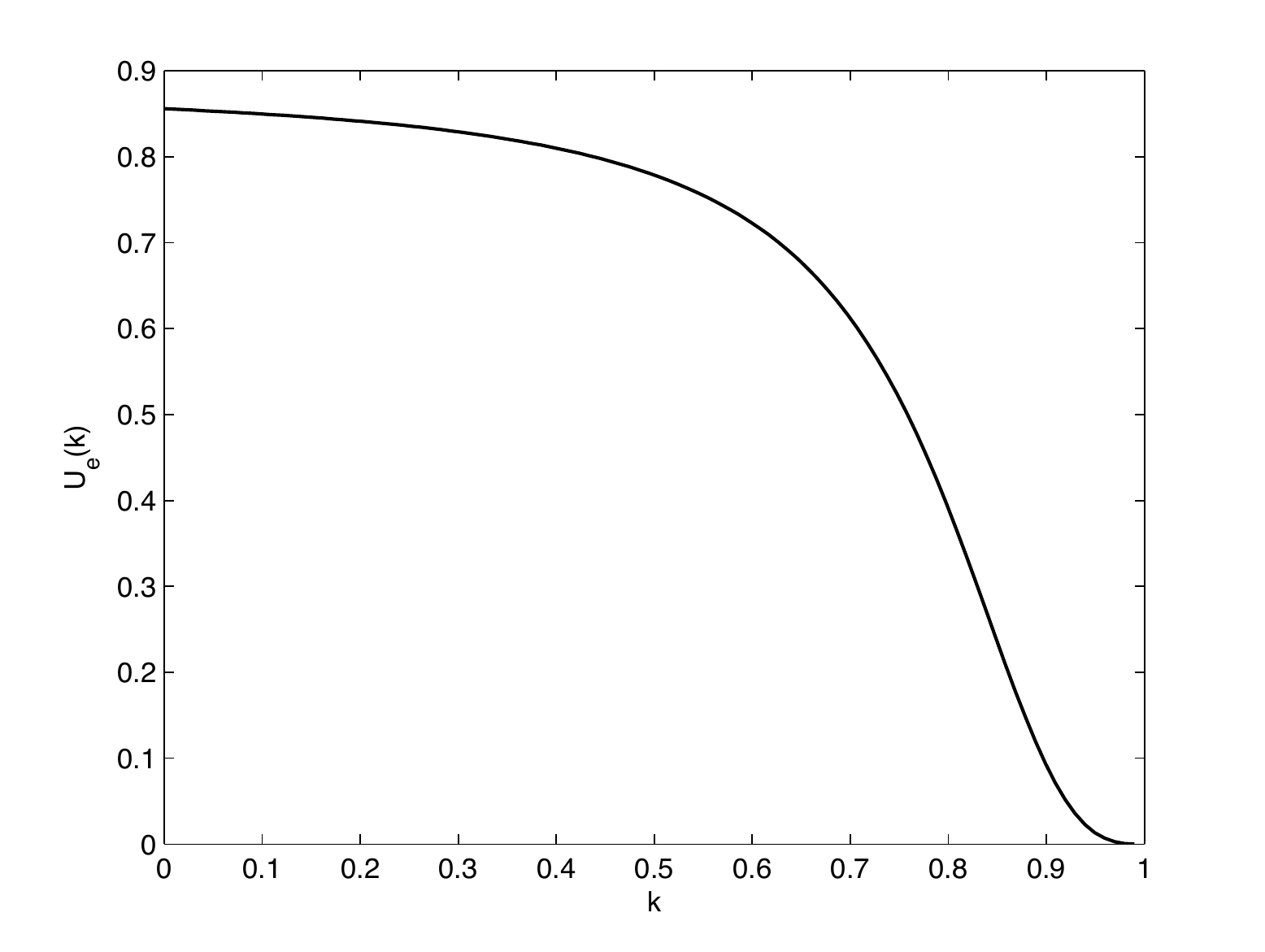,height=.3\textheight,width=.7\textwidth}
\end{center}
\caption{Function $u^e(k)$}
\label{uek}
\end{figure}

Using $P_B$ and $\nu$ described above we compute for fixed $k$ the unique stationary solution
of the homogeneous kinetic equation and
 the function $u^e(k)$
following Section 
\ref{Stationary Distributions and multi-valued Fundamental Diagrams}.
The dependence of  $u^e$ on $k$ is plotted in Figure
\ref{uek}.

\subsection{The multi-valued fundamental diagrams}

\begin{figure}
\begin{center}
\epsfig{file=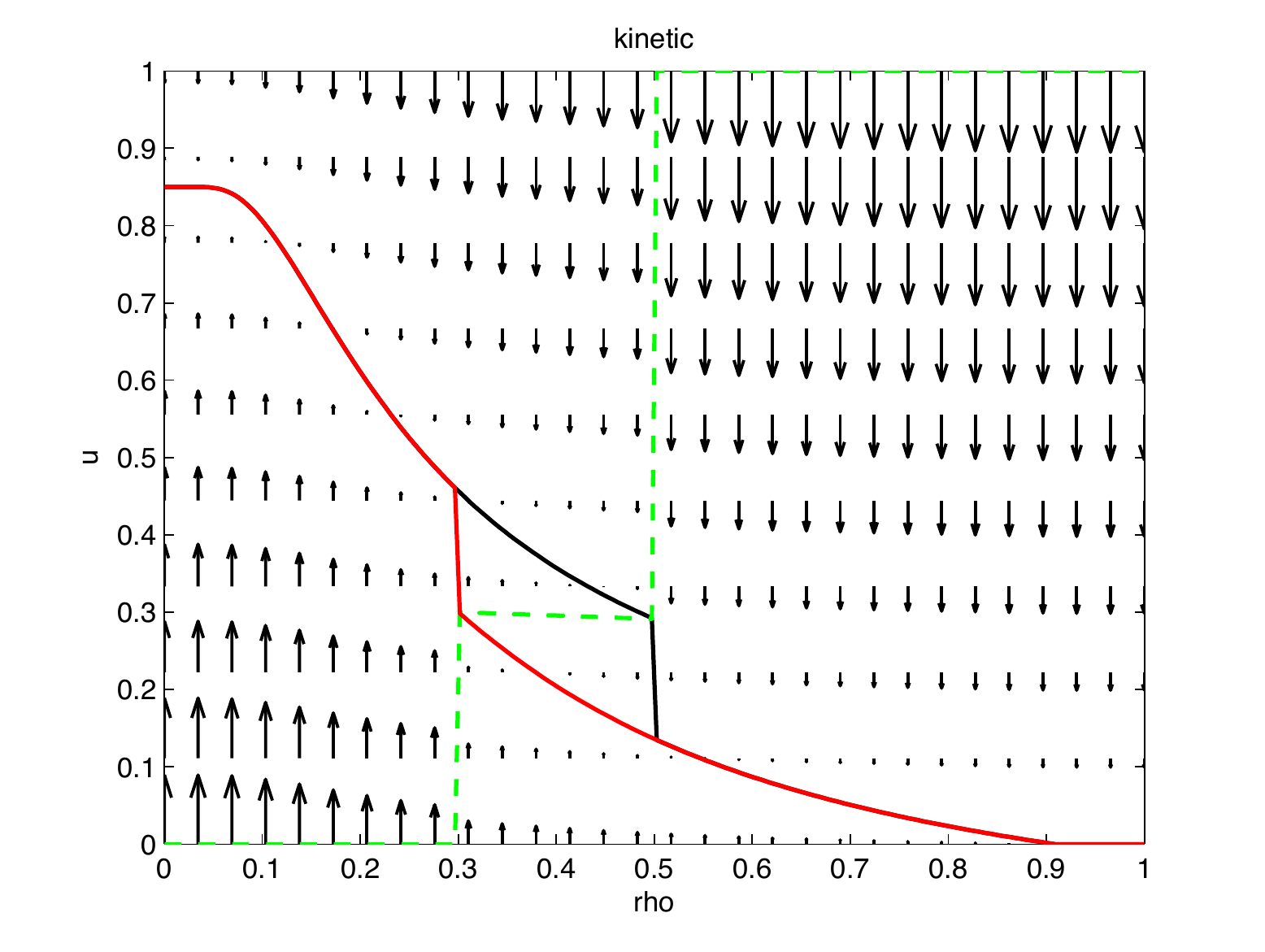,height=.25\textheight,width=.45\textwidth}
\epsfig{file=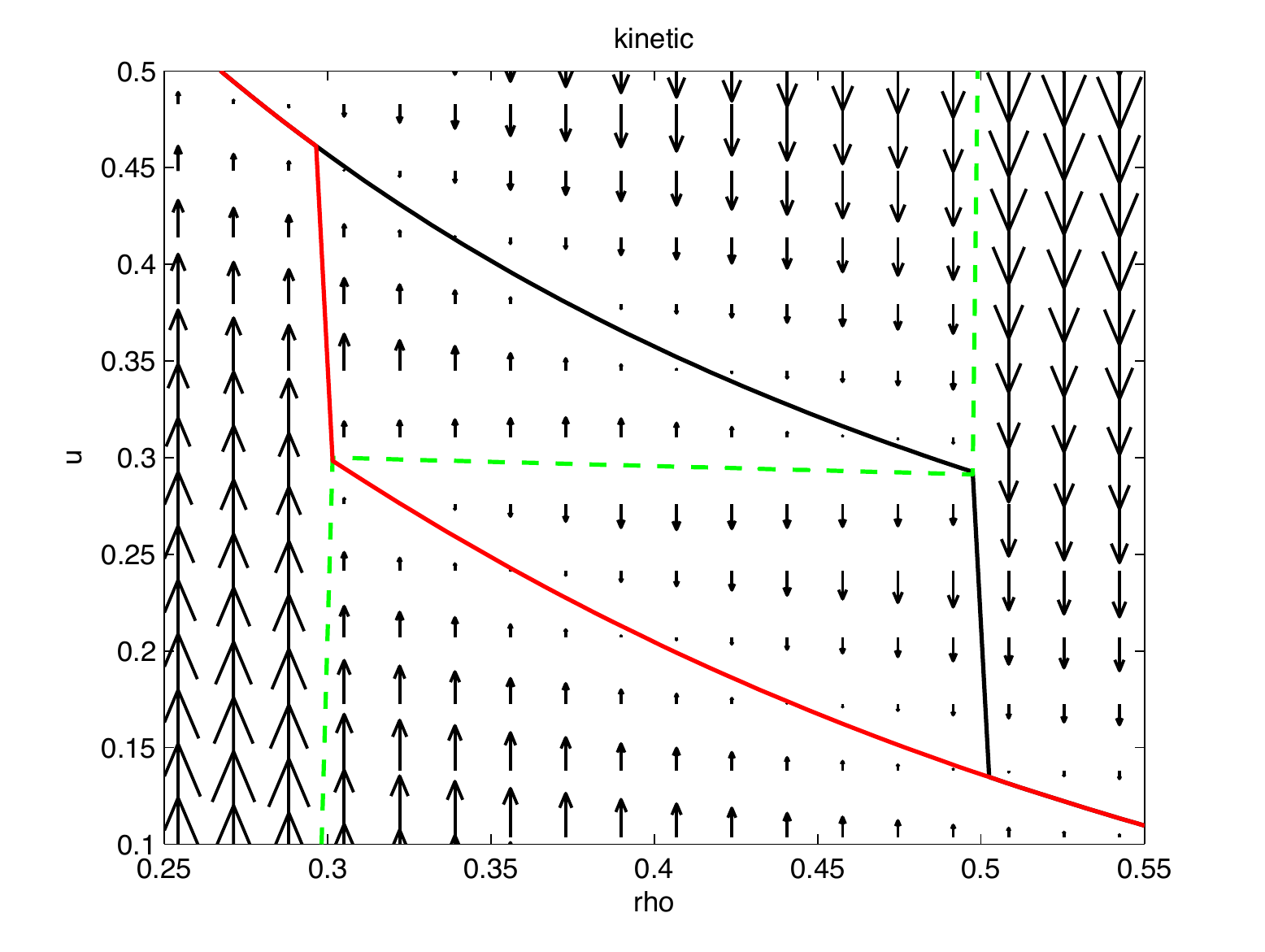,height=.25\textheight,width=.45\textwidth}
\end{center}
\caption{ $U(\rho,u)-u$ for the kinetic model. On the Right: a zoom of the multi-valued region.}
\label{uekinetic}
\end{figure}

\begin{figure}
\begin{center}
\epsfig{file=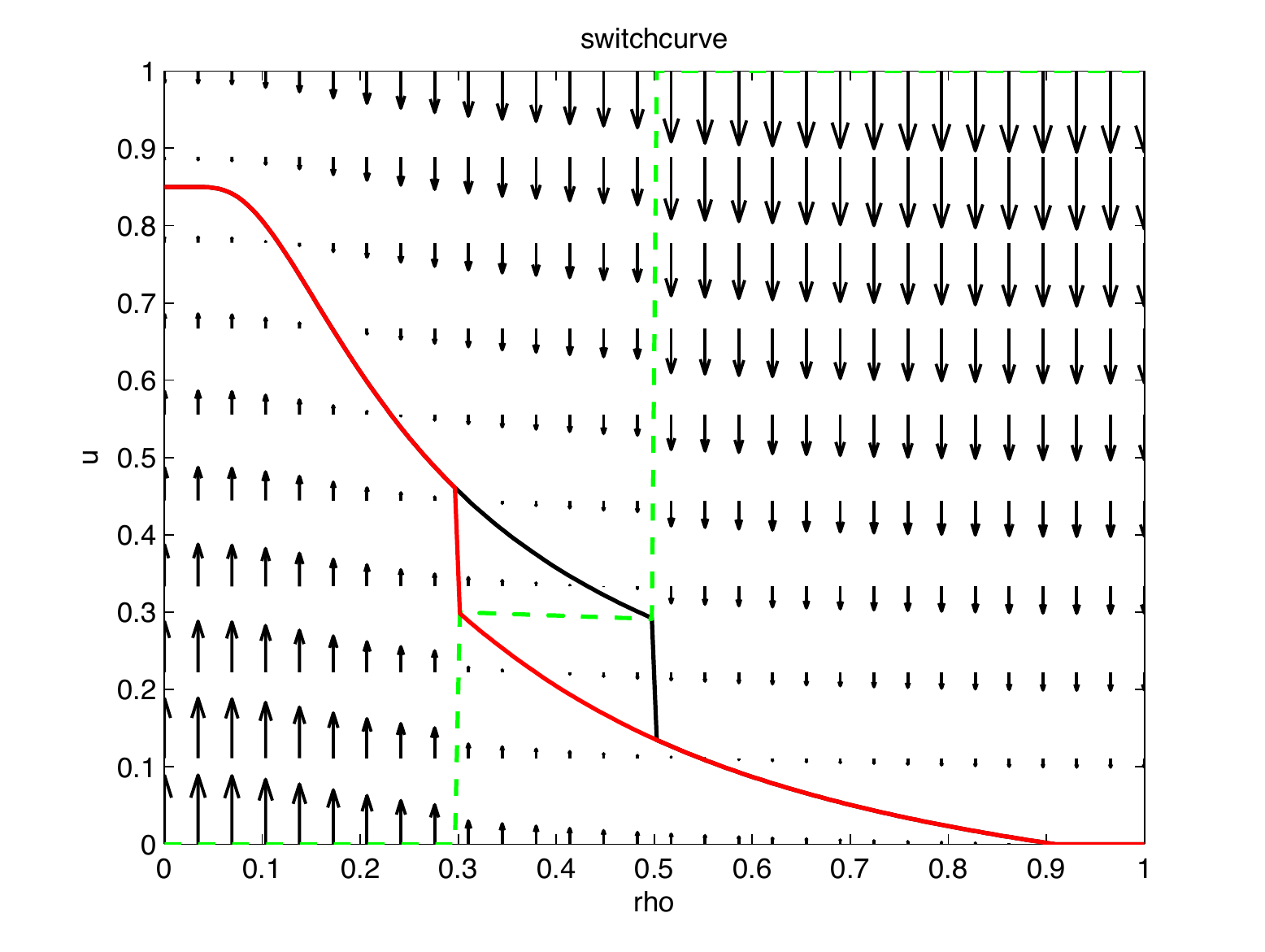,height=.25\textheight,width=.45\textwidth}
\epsfig{file=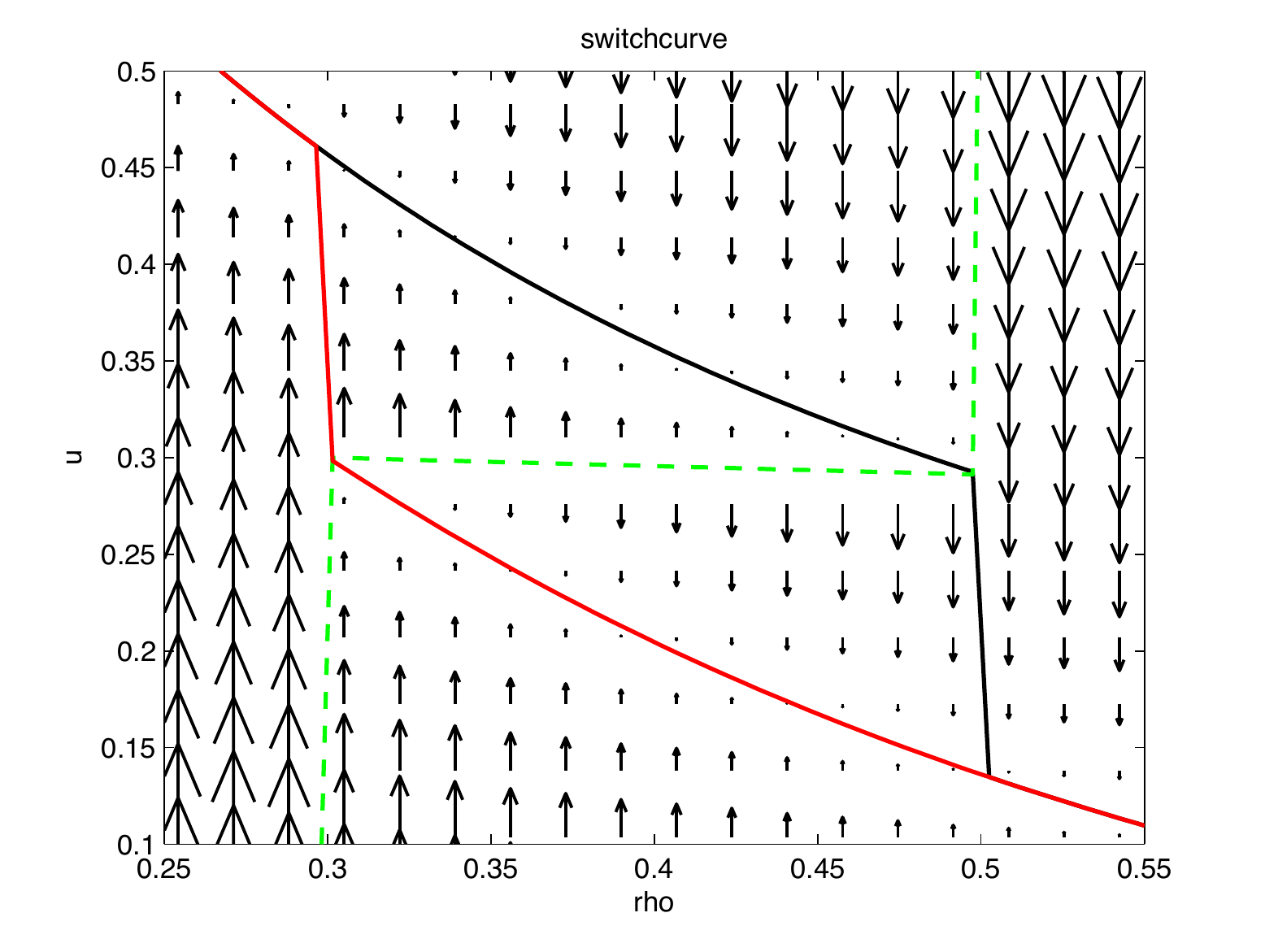,height=.25\textheight,width=.45\textwidth}
\end{center}
\caption{ $U(\rho,u)-u$ for the switching curve  model. On the Right: a zoom of the multi-valued region.}
\label{ueswitch}
\end{figure}

\begin{figure}
\begin{center}
\epsfig{file=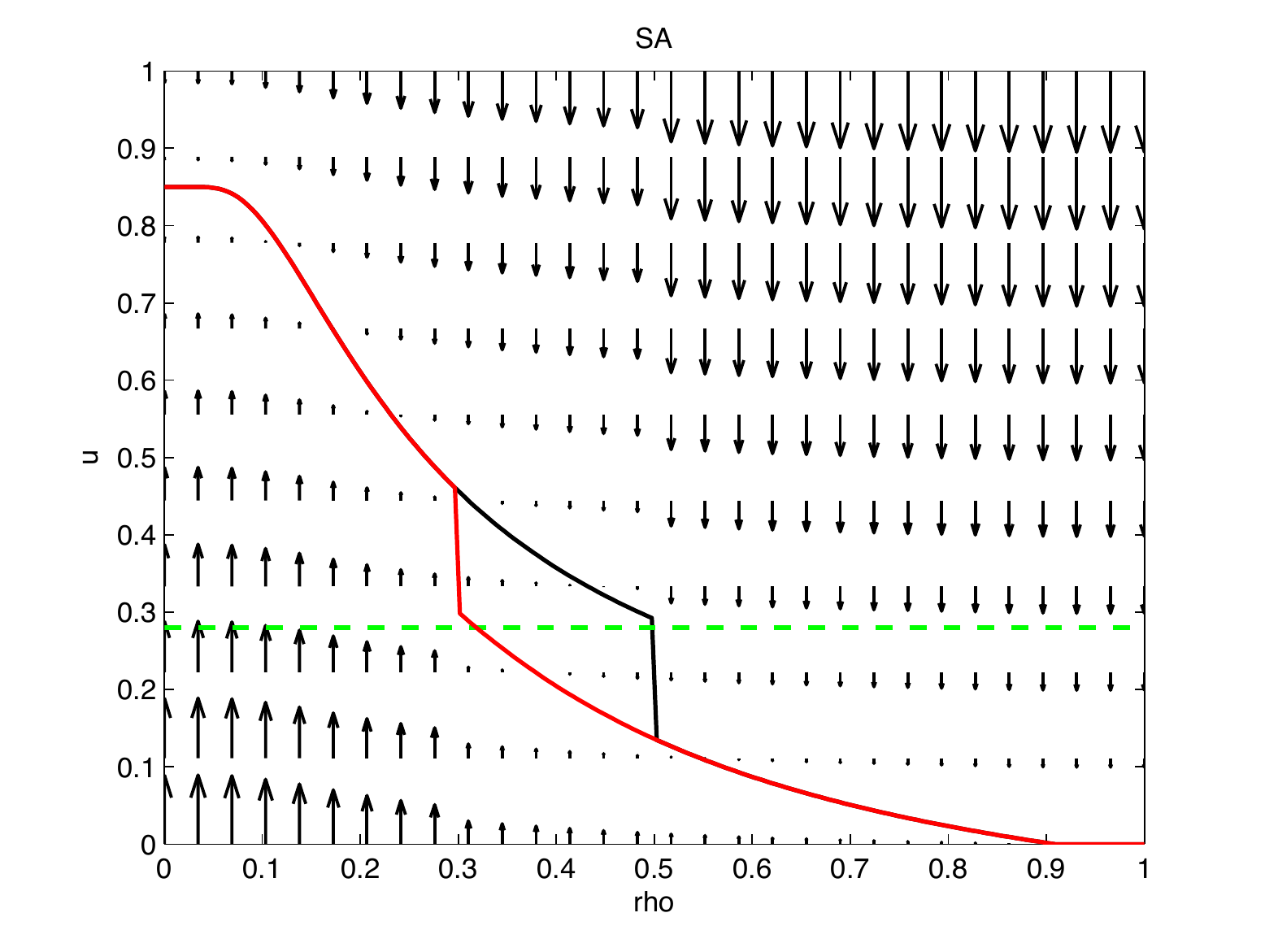,height=.25\textheight,width=.45\textwidth}
\epsfig{file=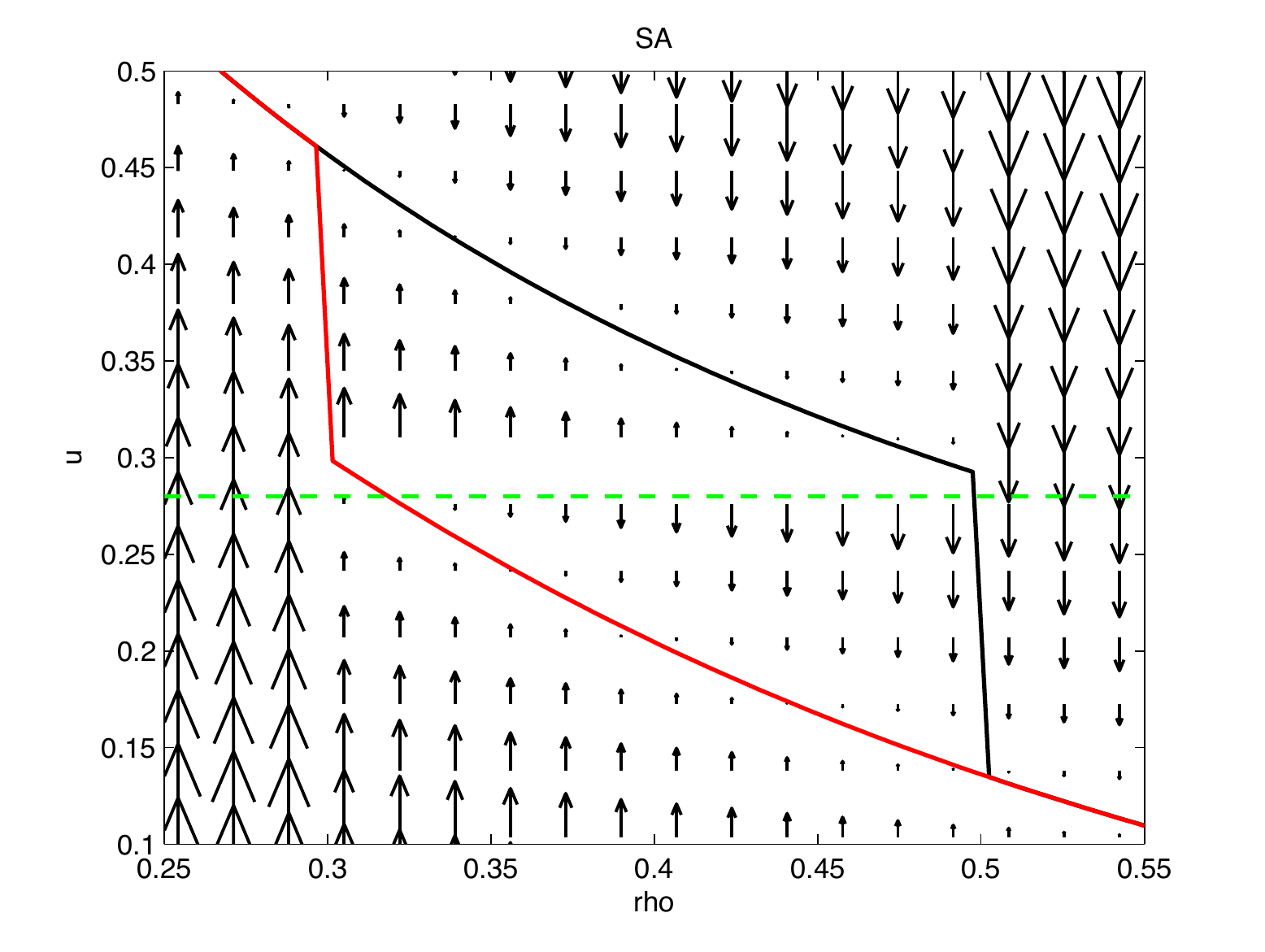,height=.25\textheight,width=.45\textwidth}
\end{center}
\caption{ $U(\rho,u)-u$ for the SA model. On the Right: a zoom of the multi-valued region.}
\label{ueSA}
\end{figure}

\begin{figure}
 \begin{center}
\epsfig{file=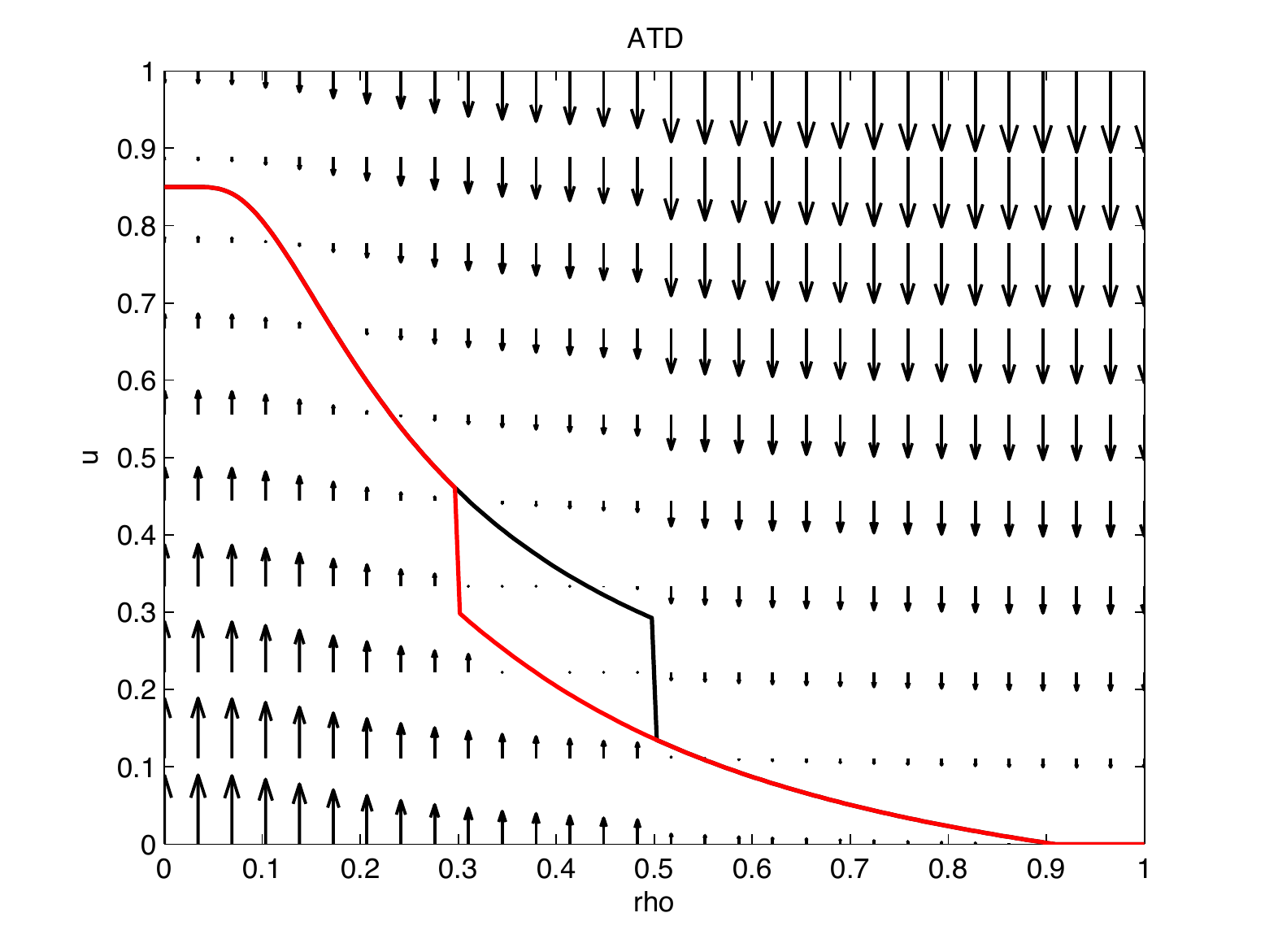,height=.25\textheight,width=.45\textwidth}
\epsfig{file=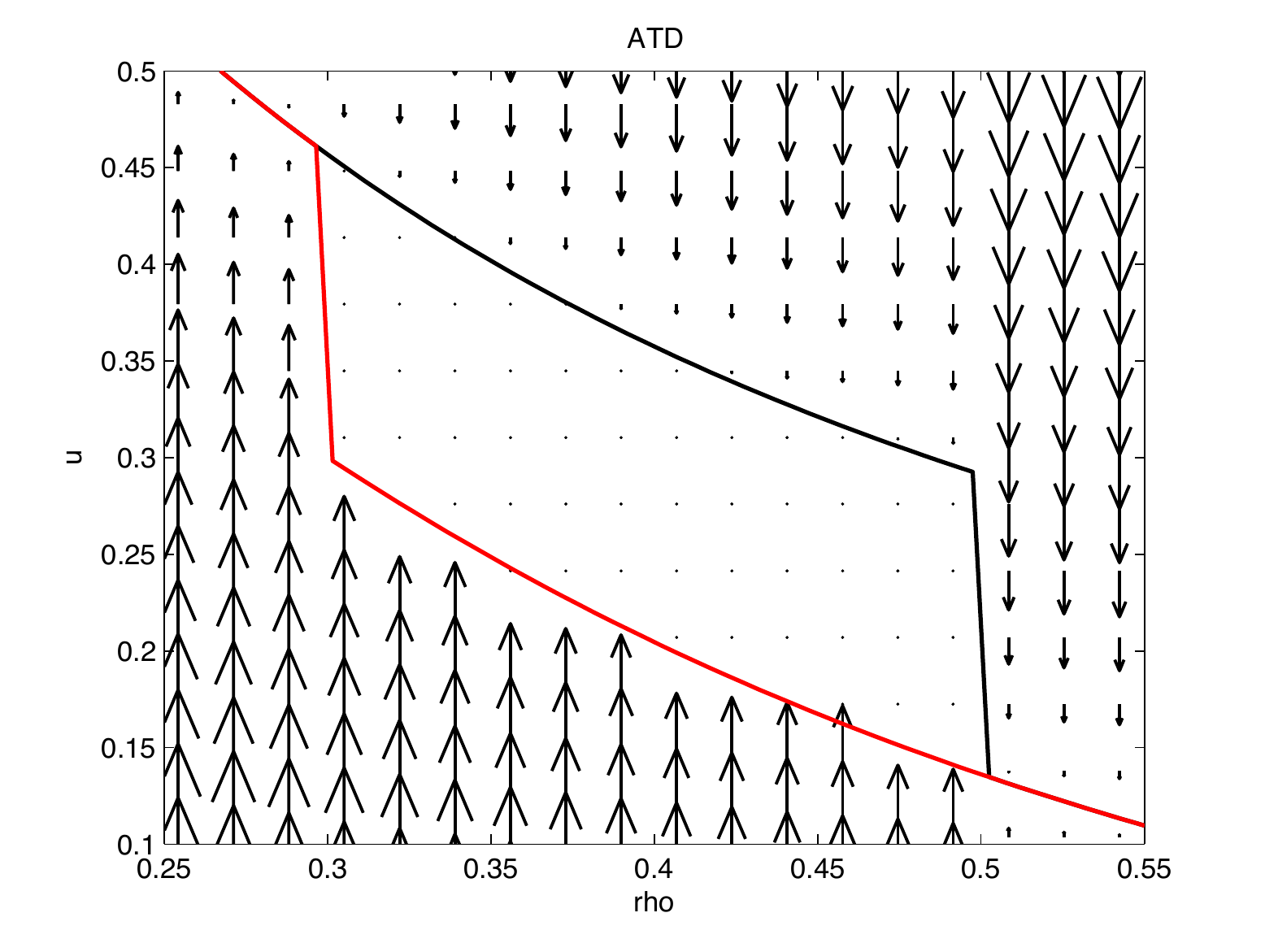,height=.25\textheight,width=.45\textwidth}
 \end{center}
 \caption{ $U(\rho,u)-u$ for the ATD model. On the Right: a zoom of the multi-valued region.}
 \label{ueATD}
 \end{figure}

In this subsection we plot the multi-valued fundamental diagrams for the four cases discussed in section \ref{Comparison of multi-phase hydrodynamic models}. 
The functions $U_1(\rho)$ and $U_2(\rho)$ for ATD and switching curve models were chosen as
\begin{eqnarray*}
U_1(\rho)\;=  U_0\mbox{tanh}\left( \frac{C_U}{T_0 U_0} \left(\frac{1}{\rho} -0.05 \right)\right)
\end{eqnarray*}
\begin{eqnarray*}
U_2(\rho)\;=  U_0^*\mbox{tanh}\left( \frac{C_U}{T_0 U_0^*} \left(\frac{1}{\rho} -1.1 \right)\right)
\end{eqnarray*}
with $U_0=0.85,\;C_U=0.45,\;U_0^*=0.5,\;T_0=2.9$. Moreover, $T=5,\;C=0.3,\;U_{sync}=0.28$ and $S=S(\rho) $ is given by a linear function connecting $U_1(\rho_f)$ and $U_2(\rho_j)$. The solutions of the  nonlinear equation
$u = U (\rho,u)$ are plotted together with a plot of the quantity $U (\rho,u) -u $  as arrows with direction $\pm u$.
Figure \ref{uekinetic}, \ref{ueswitch},  \ref{ueSA}  and \ref{ueATD} show the multi-valued fundamental diagram
(speed-density relation) for the different models. In each figure, a zoom of the multi-valued region is shown.

In all cases the values for $\rho_f$ and $\rho_j$ are chosen as $0.3$ and $0.5$ respectively such that for $\rho < \rho_f$ we have only one steady solution.  For $\rho_f < \rho < \rho_j$ three (infinitely many for ATD) solutions
 exist.  And for the region $\rho > \rho_j$ again only one solution exists.

\begin{figure}
\begin{center}
\epsfig{file=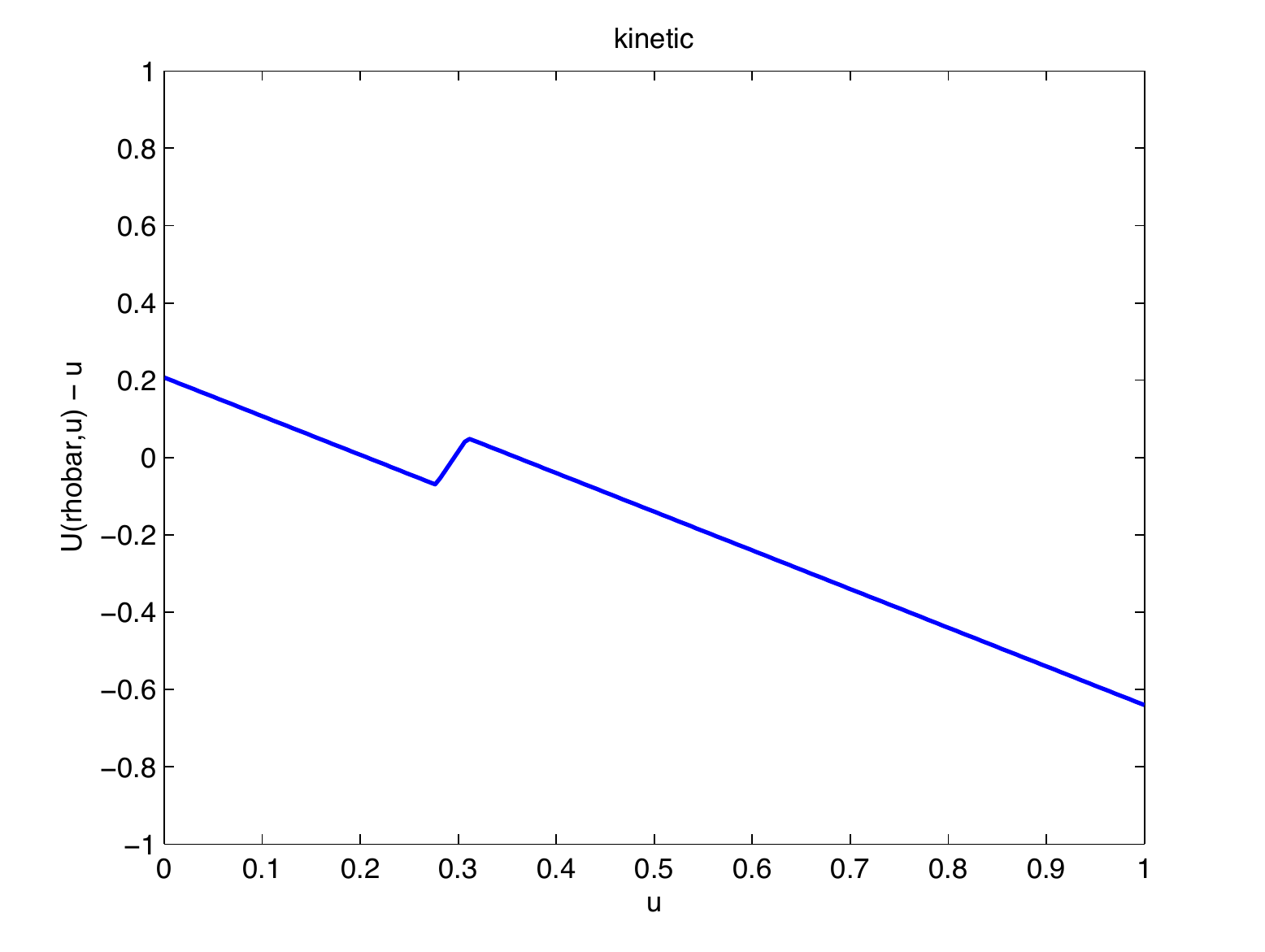,height=.15\textheight,width=.35\textwidth}
\epsfig{file=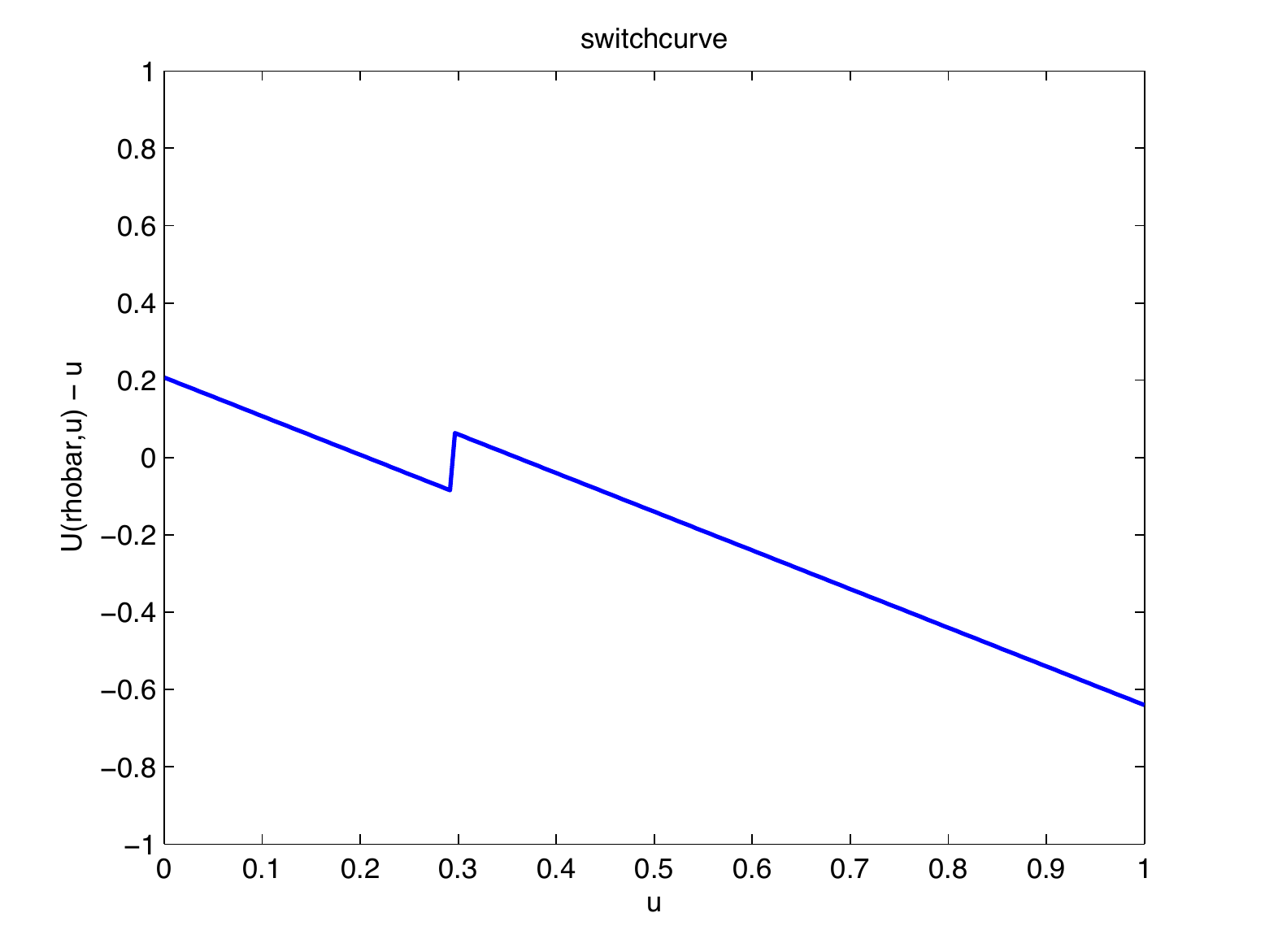,height=.15\textheight,width=.35\textwidth}
\epsfig{file=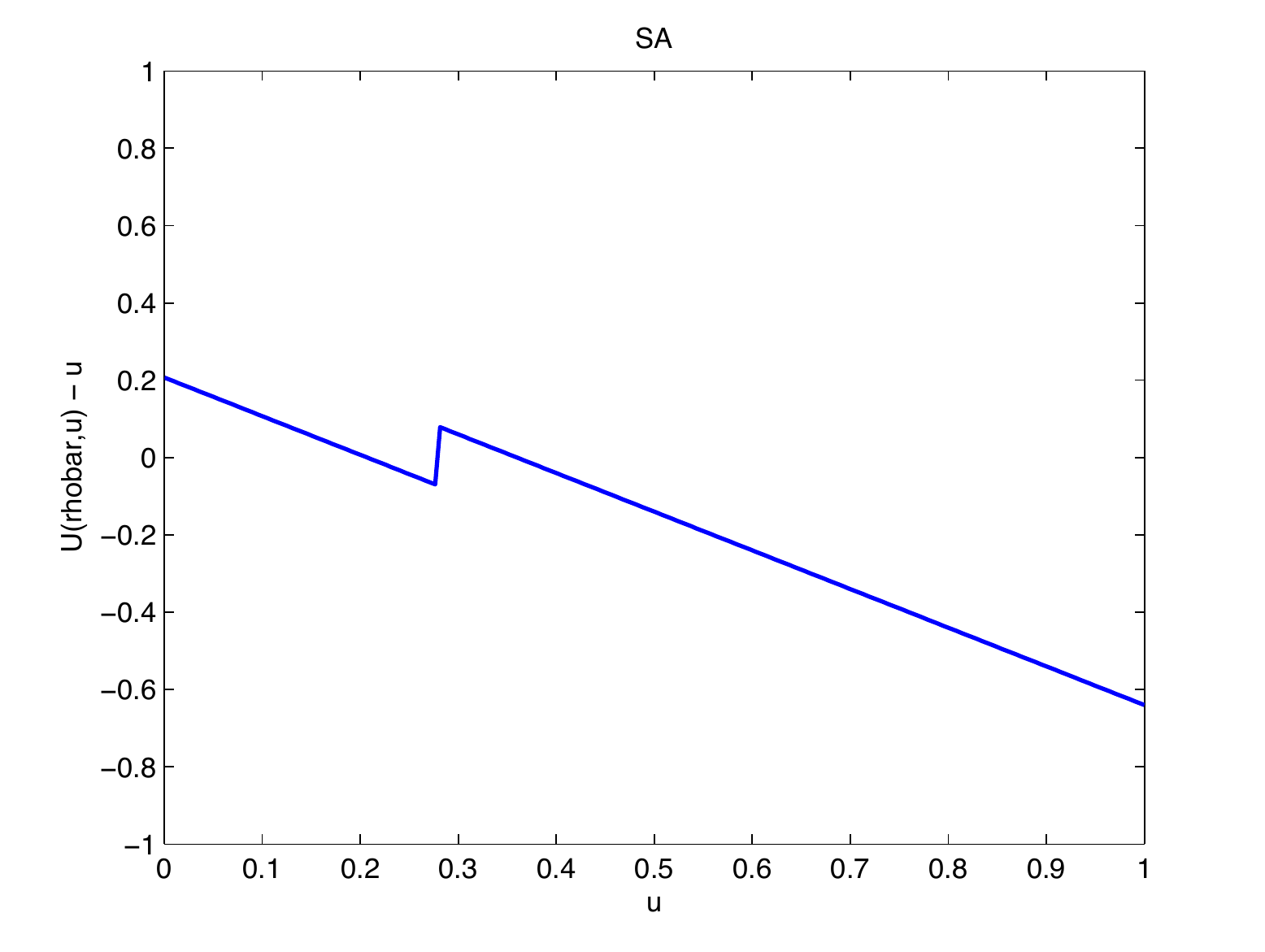,height=.15\textheight,width=.35\textwidth}
\epsfig{file=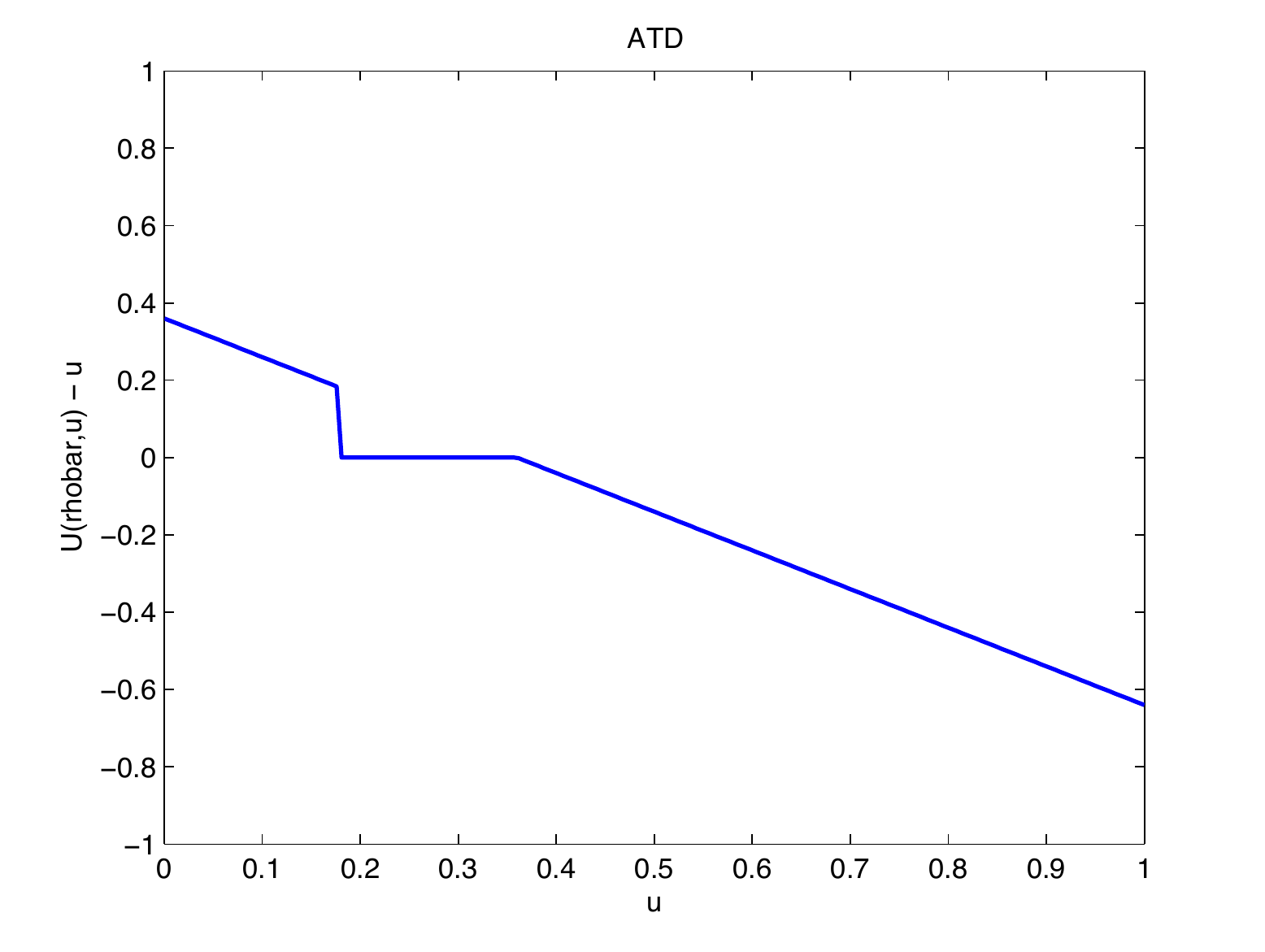,height=.15\textheight,width=.35\textwidth}
\end{center}
\caption{$U(\bar \rho,u)-u$ for the kinetic, switching curve, SA and ATD models.}
\label{uekineticcut}
\end{figure}

Moreover,   Figure \ref{uekineticcut} shows the values of $U(\rho,u)-u$ for a fixed value $\rho = \bar \rho$  with $\rho_f < \bar \rho < \rho_j$.

\subsection{Numerical solution of multi-phase macroscopic equations}

\begin{figure}
\begin{center}
 \epsfig{file=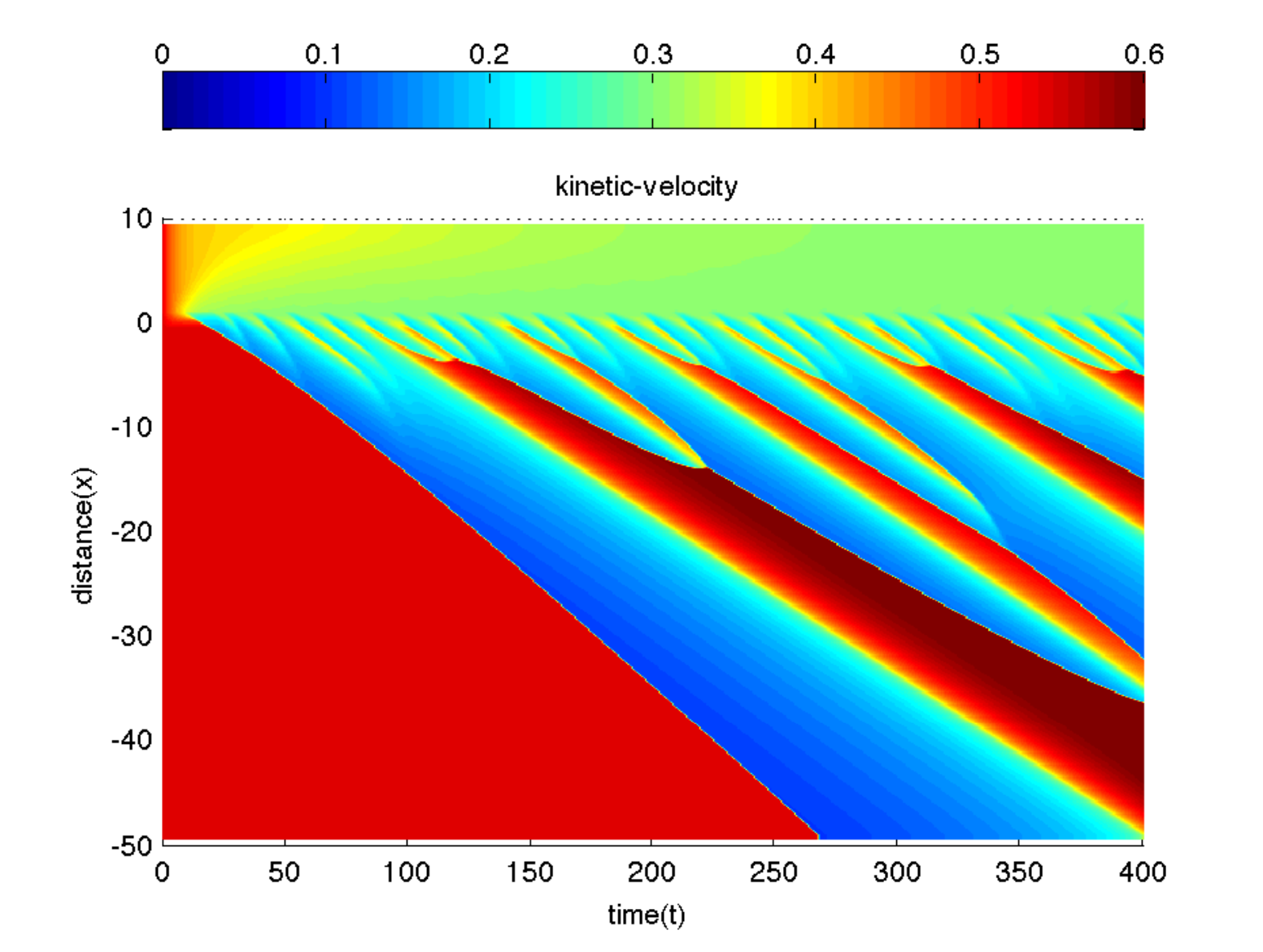,height=.25\textheight,width=.45\textwidth}
 \epsfig{file=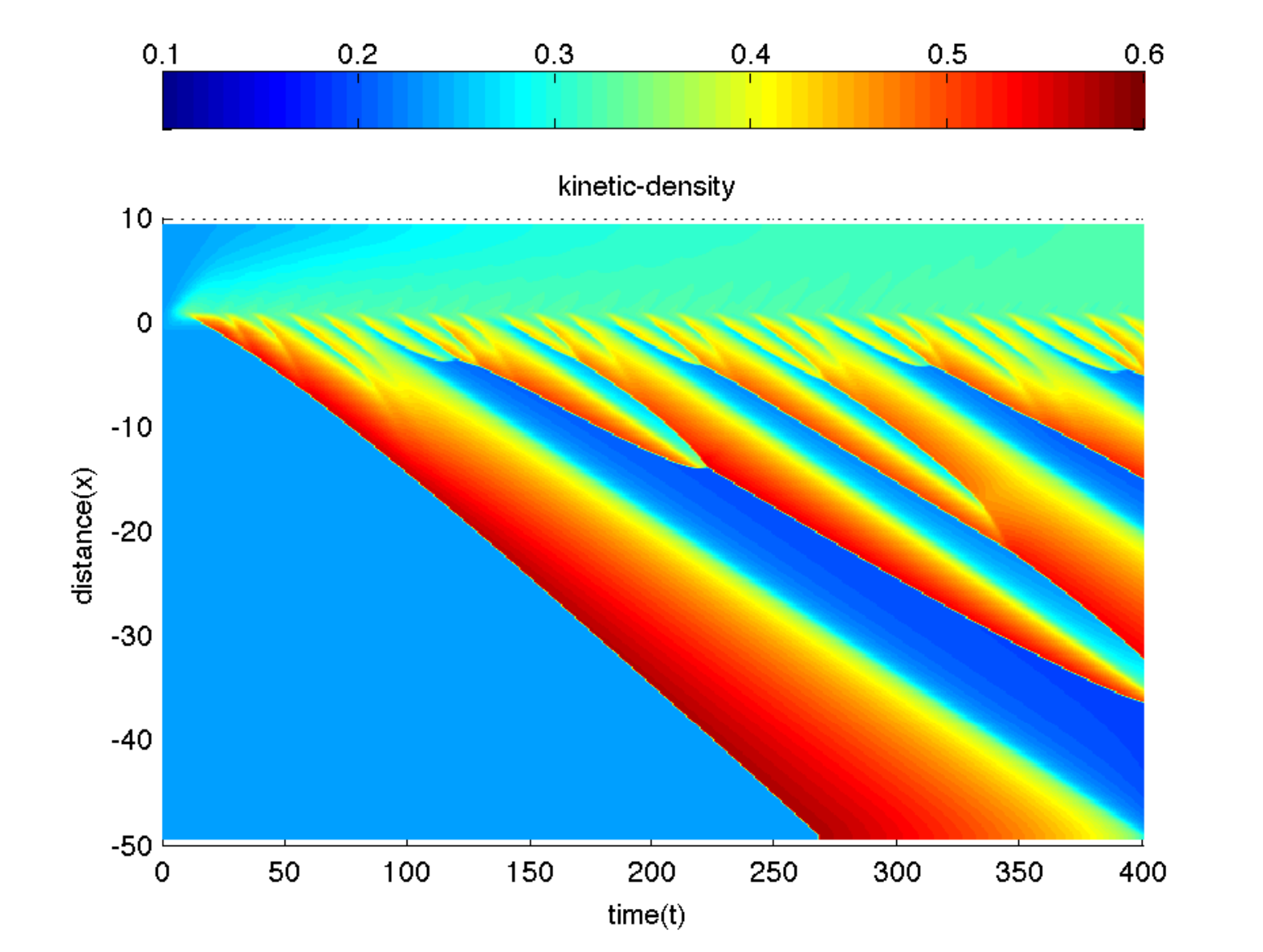,height=.25\textheight,width=.45\textwidth}
 \end{center}
\caption{Spatio-temporal congested traffic pattern - velocity(\textit{Left}) and density(\textit{Right})  for the kinetic model.}
\label{rhofluxkinetic}
\end{figure}

\begin{figure}
\begin{center}
 \epsfig{file=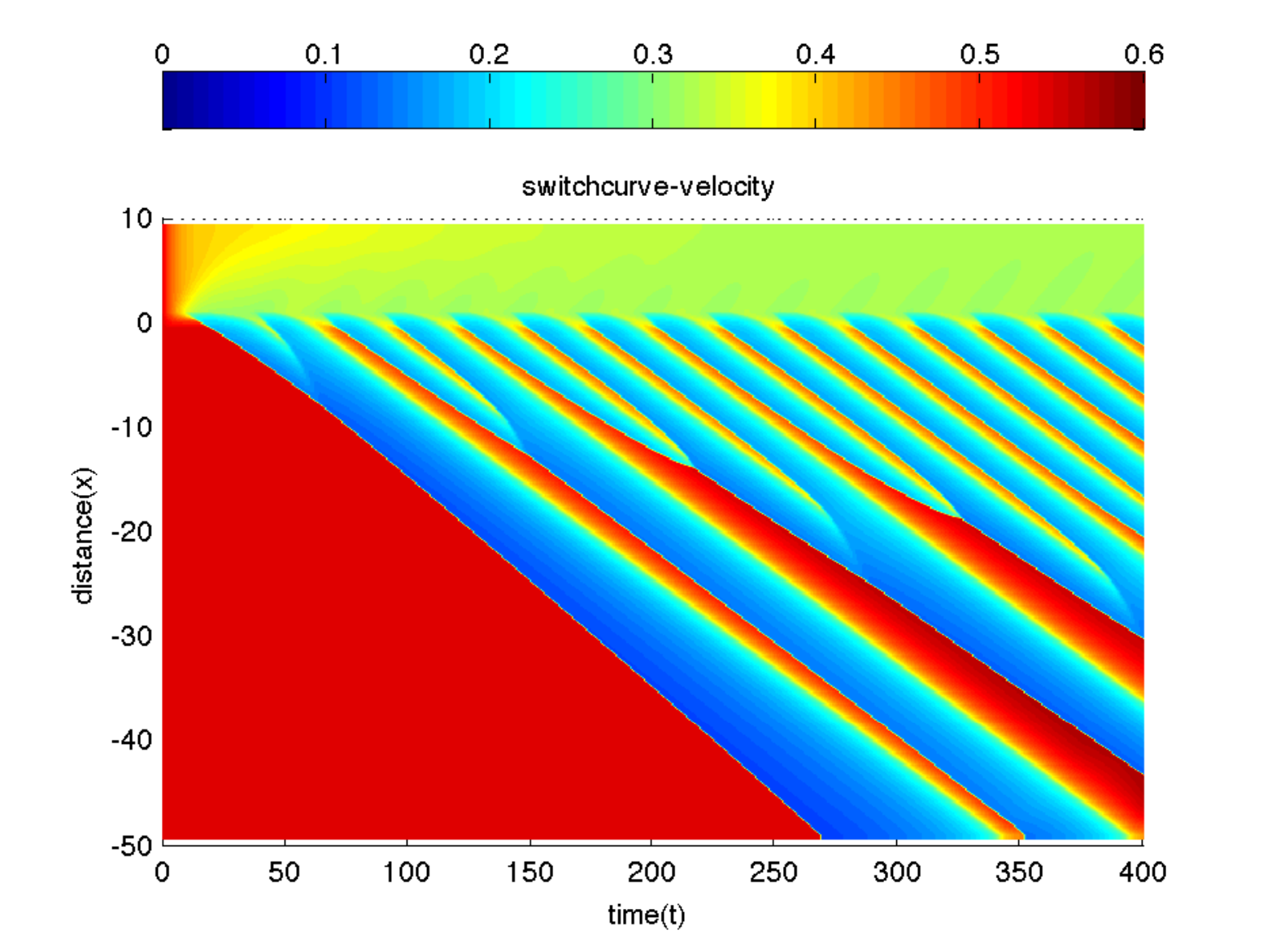,height=.25\textheight,width=.45\textwidth}
  \epsfig{file=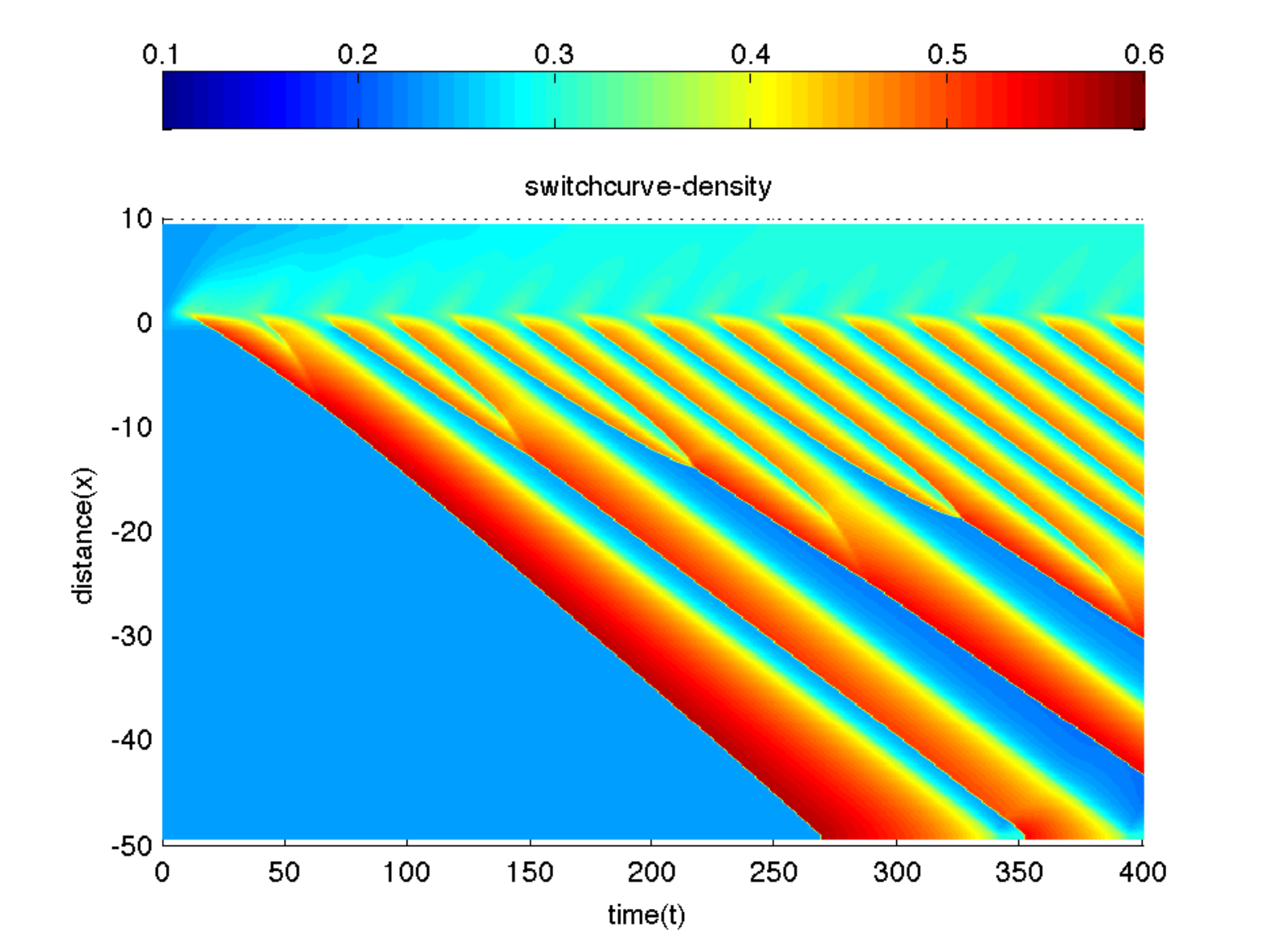,height=.25\textheight,width=.45\textwidth}
\end{center}
\caption{Spatio-temporal congested traffic pattern -  velocity(\textit{Left}) and density(\textit{Right})  for the switching curve model.}
\label{rhofluxswitch}
\end{figure}

\begin{figure}
\begin{center}
 \epsfig{file=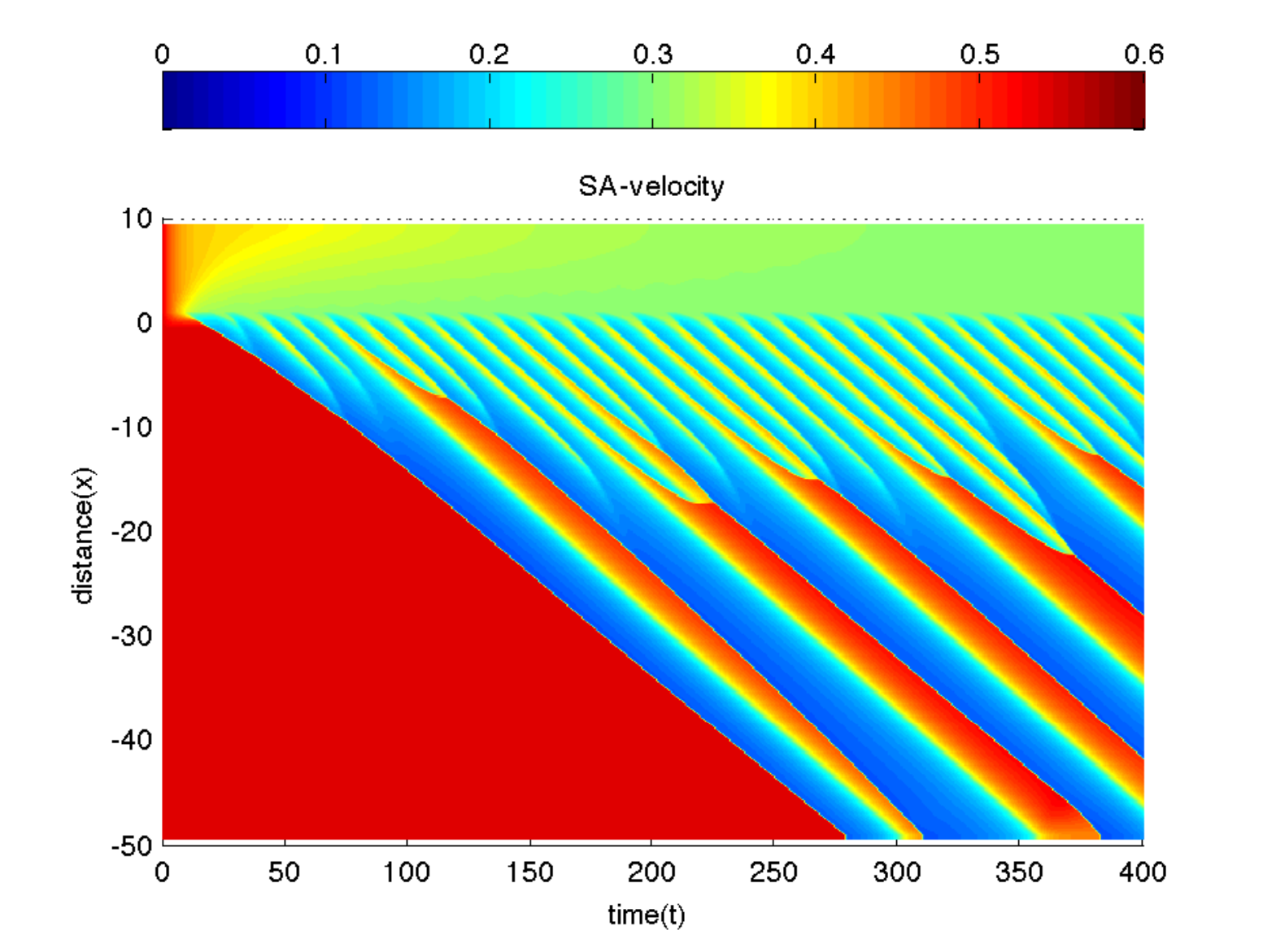,height=.25\textheight,width=.45\textwidth}
  \epsfig{file=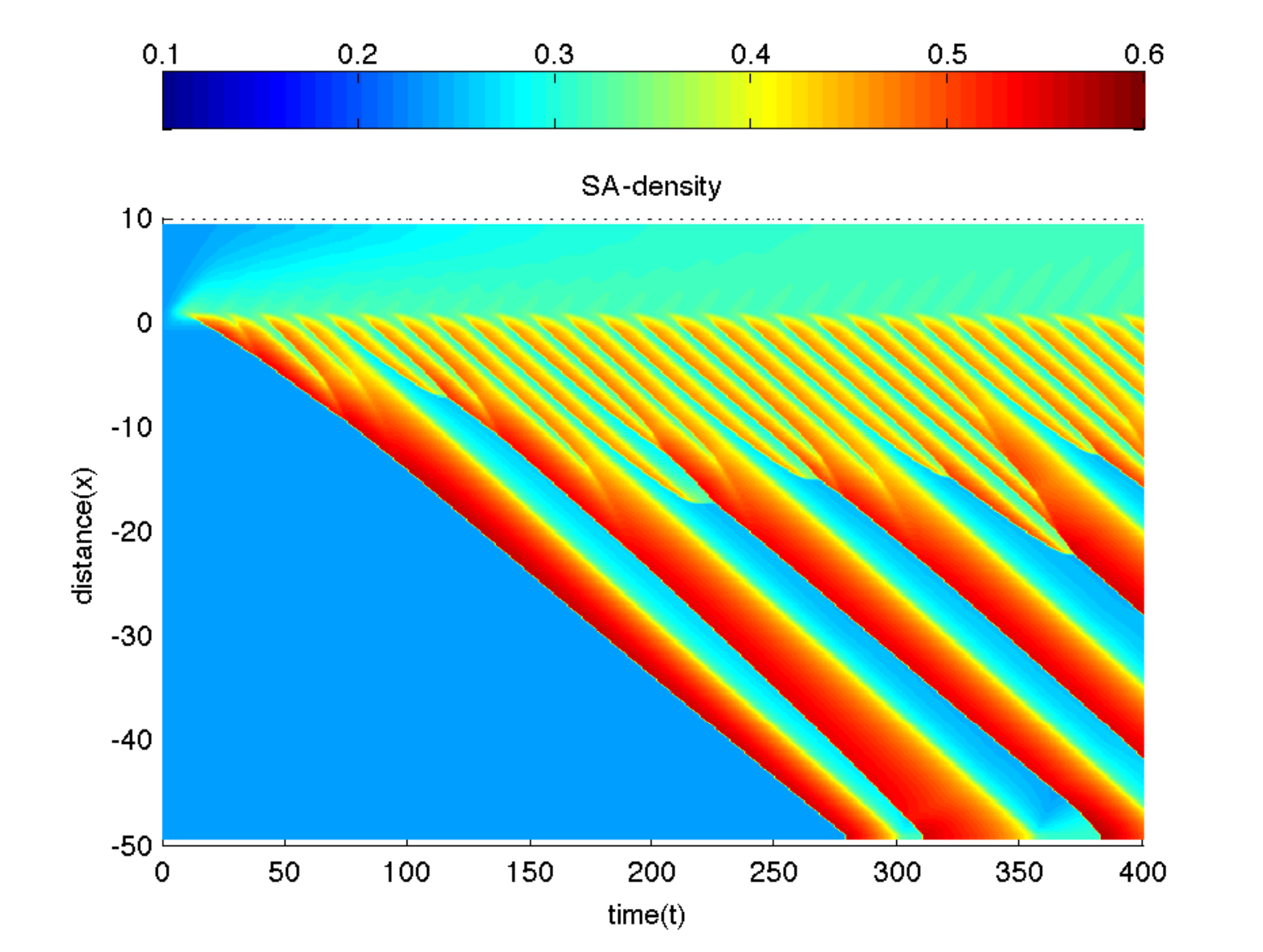,height=.25\textheight,width=.45\textwidth}
\end{center}
\caption{Spatio-temporal congested traffic pattern - velocity(\textit{Left}) and density(\textit{Right}) for the SA model.}
\label{rhofluxSA}
\end{figure}

\begin{figure}
\begin{center}
 \epsfig{file=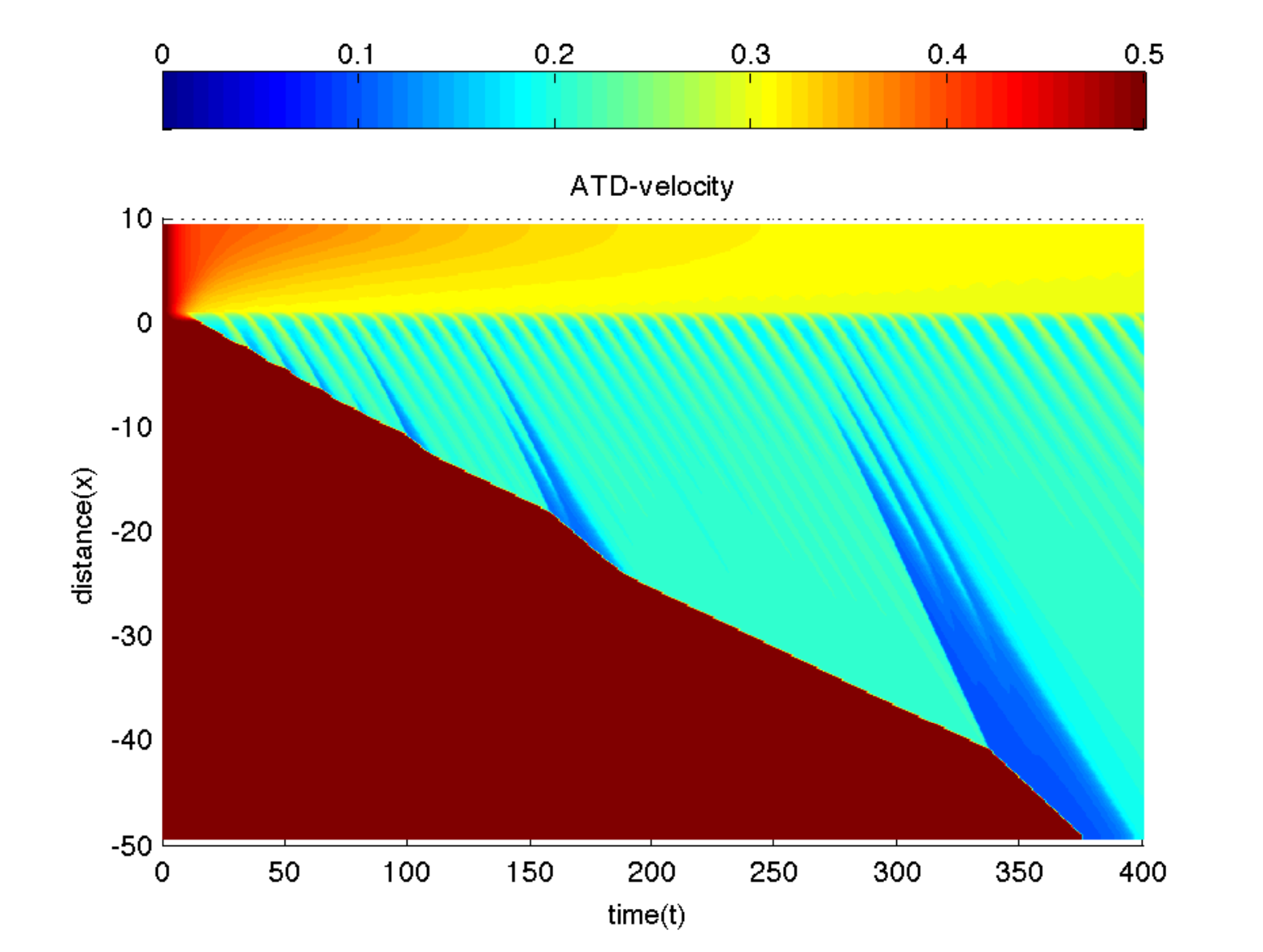,height=.25\textheight,width=.45\textwidth}
  \epsfig{file=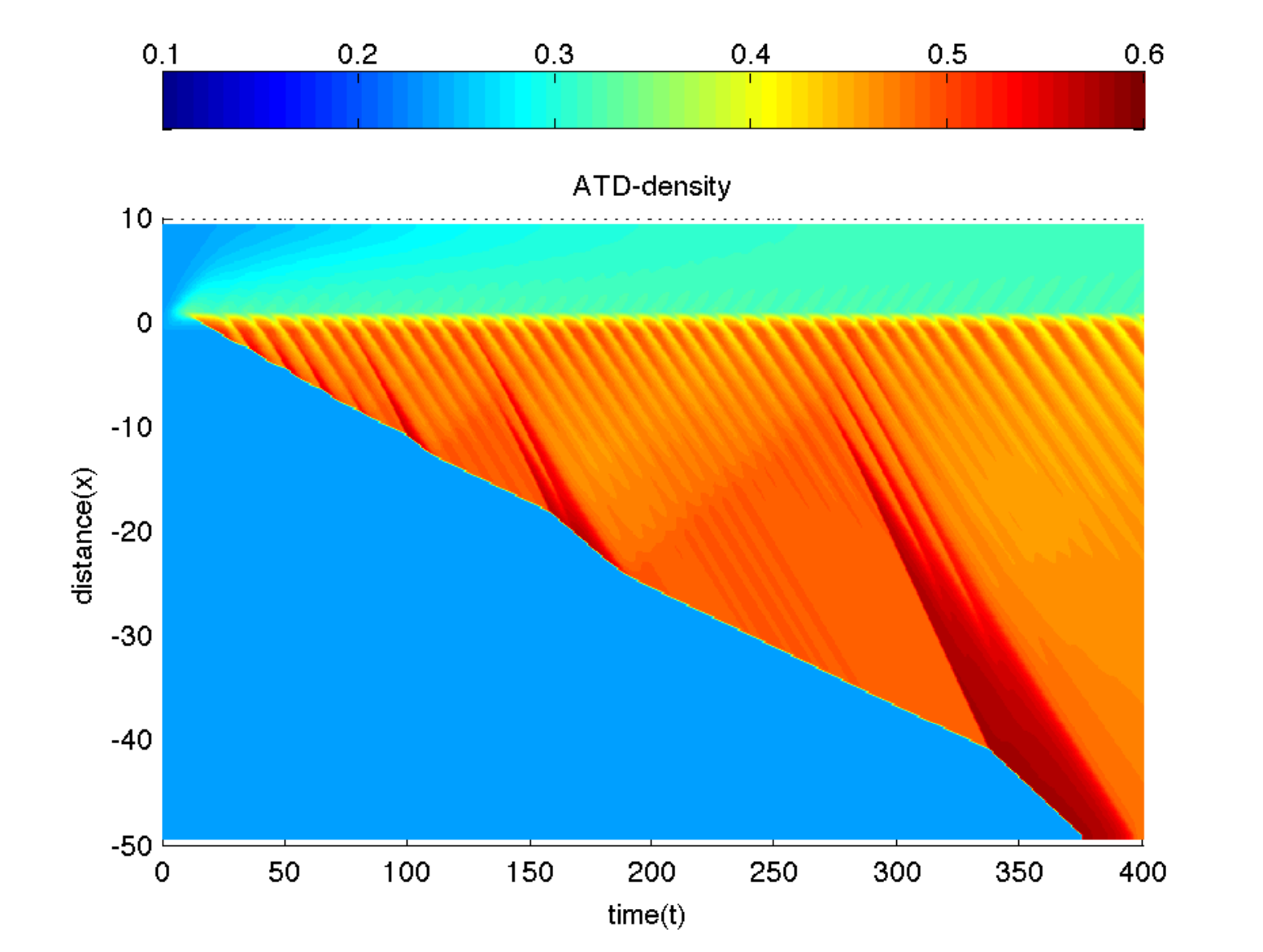,height=.25\textheight,width=.45\textwidth}
\end{center}
\caption{Spatio-temporal congesting traffic pattern - velocity(\textit{Left}) and density(\textit{Right}) for the ATD model.}
\label{rhofluxATD}
\end{figure}

\begin{figure}
\begin{center}
\epsfig{file=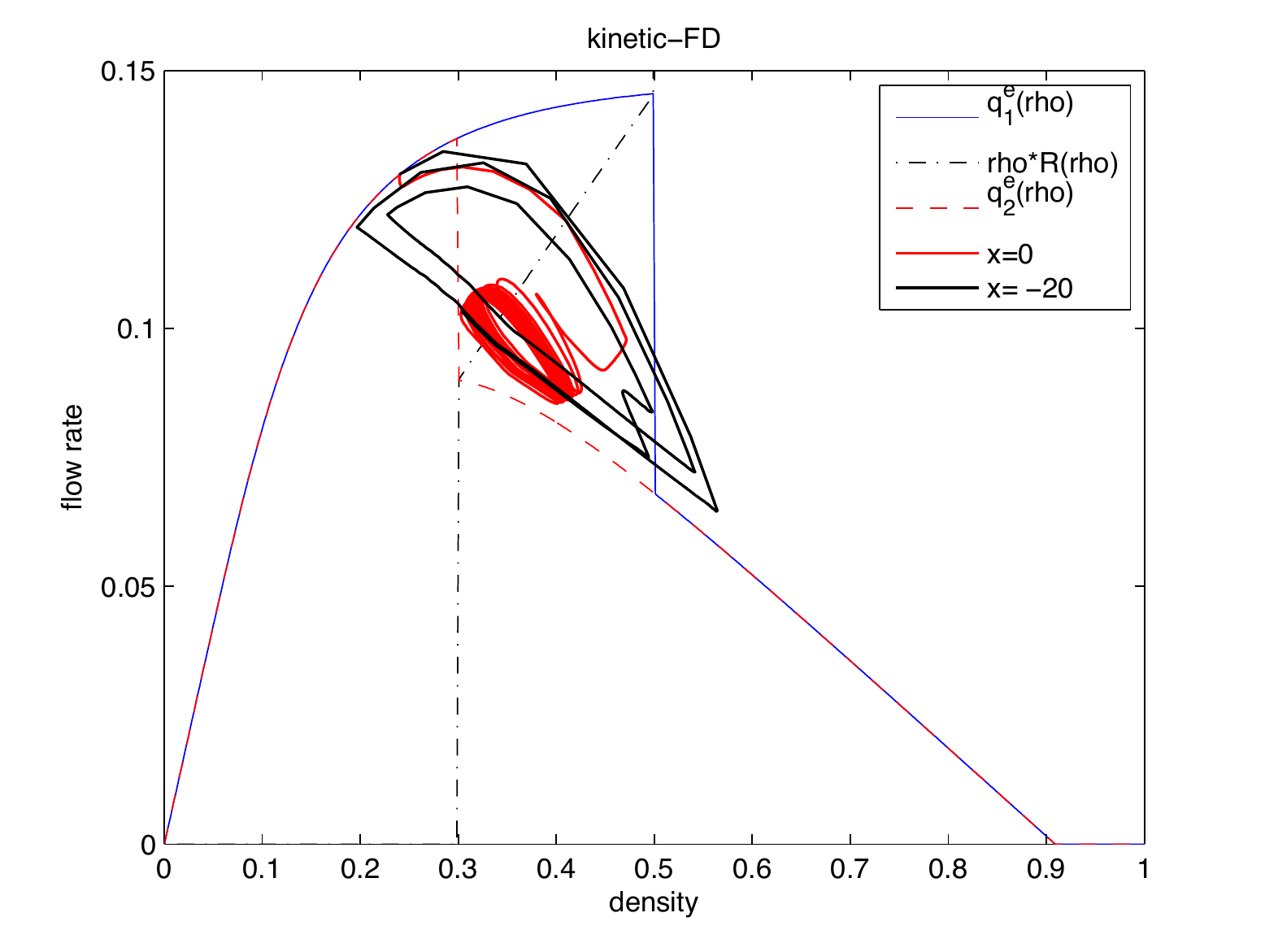,height=.25\textheight,width=.45\textwidth}
\epsfig{file=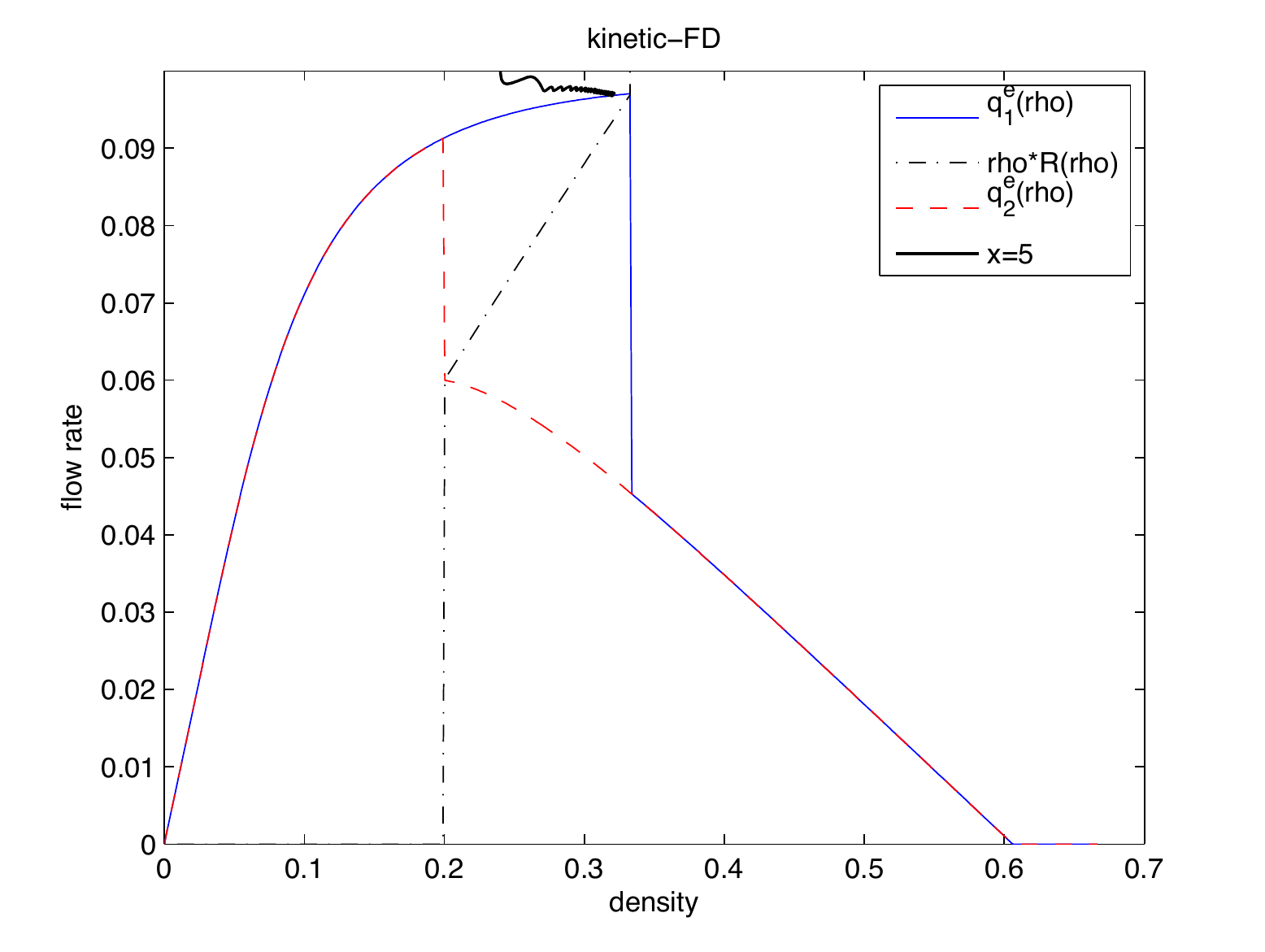,height=.25\textheight,width=.45\textwidth}
\epsfig{file=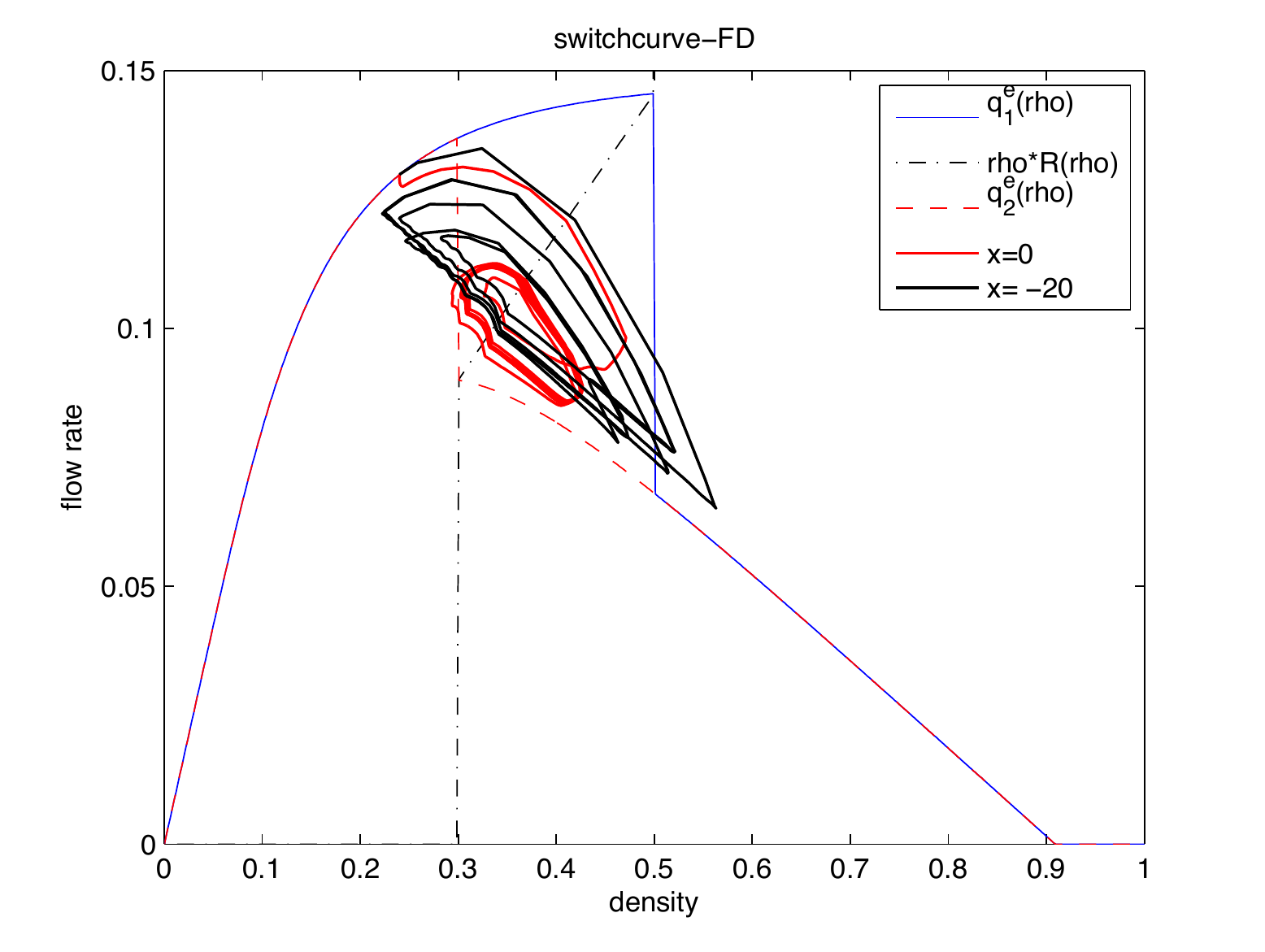,height=.25\textheight,width=.45\textwidth}
\epsfig{file=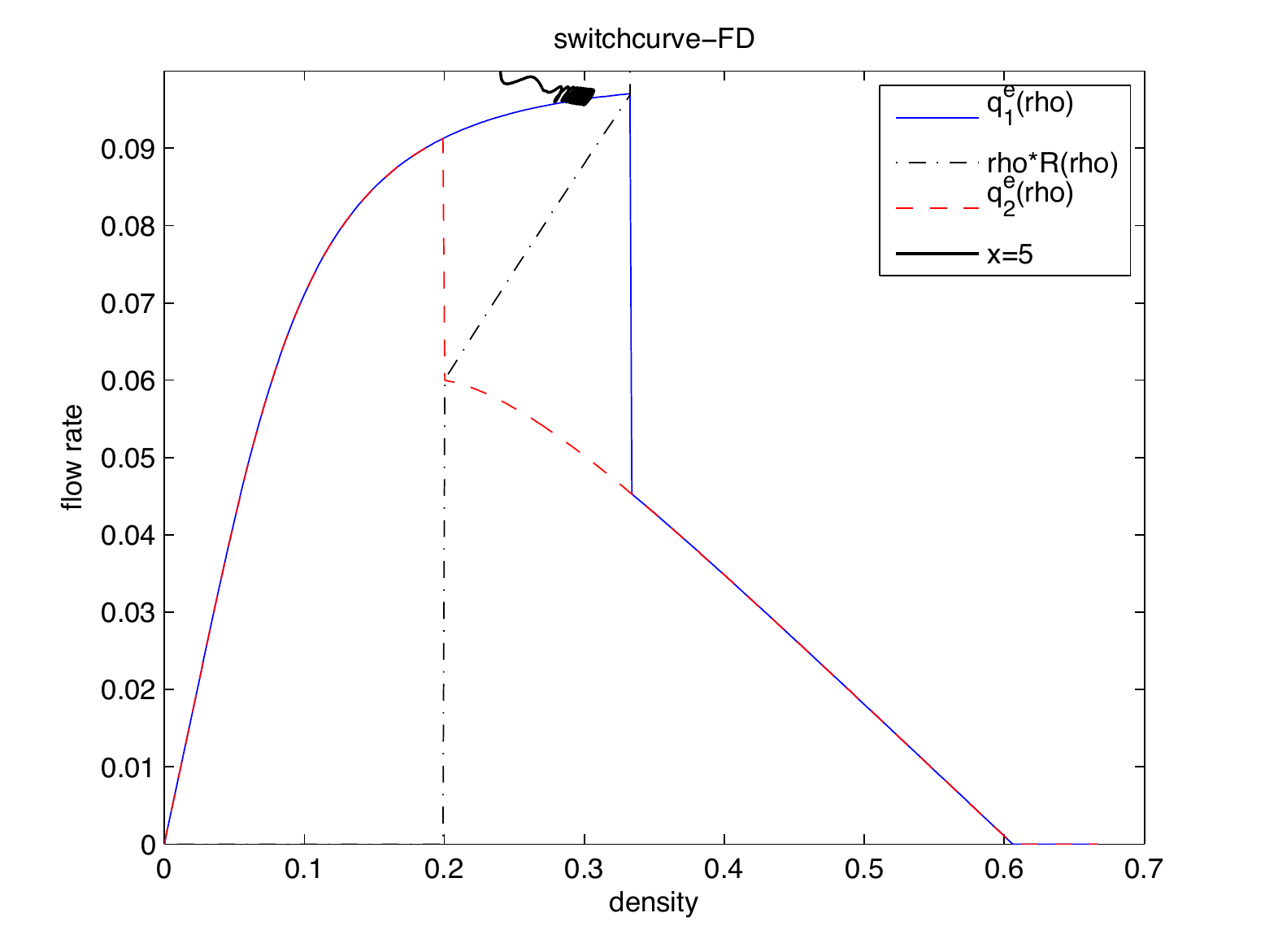,height=.25\textheight,width=.45\textwidth}
\epsfig{file=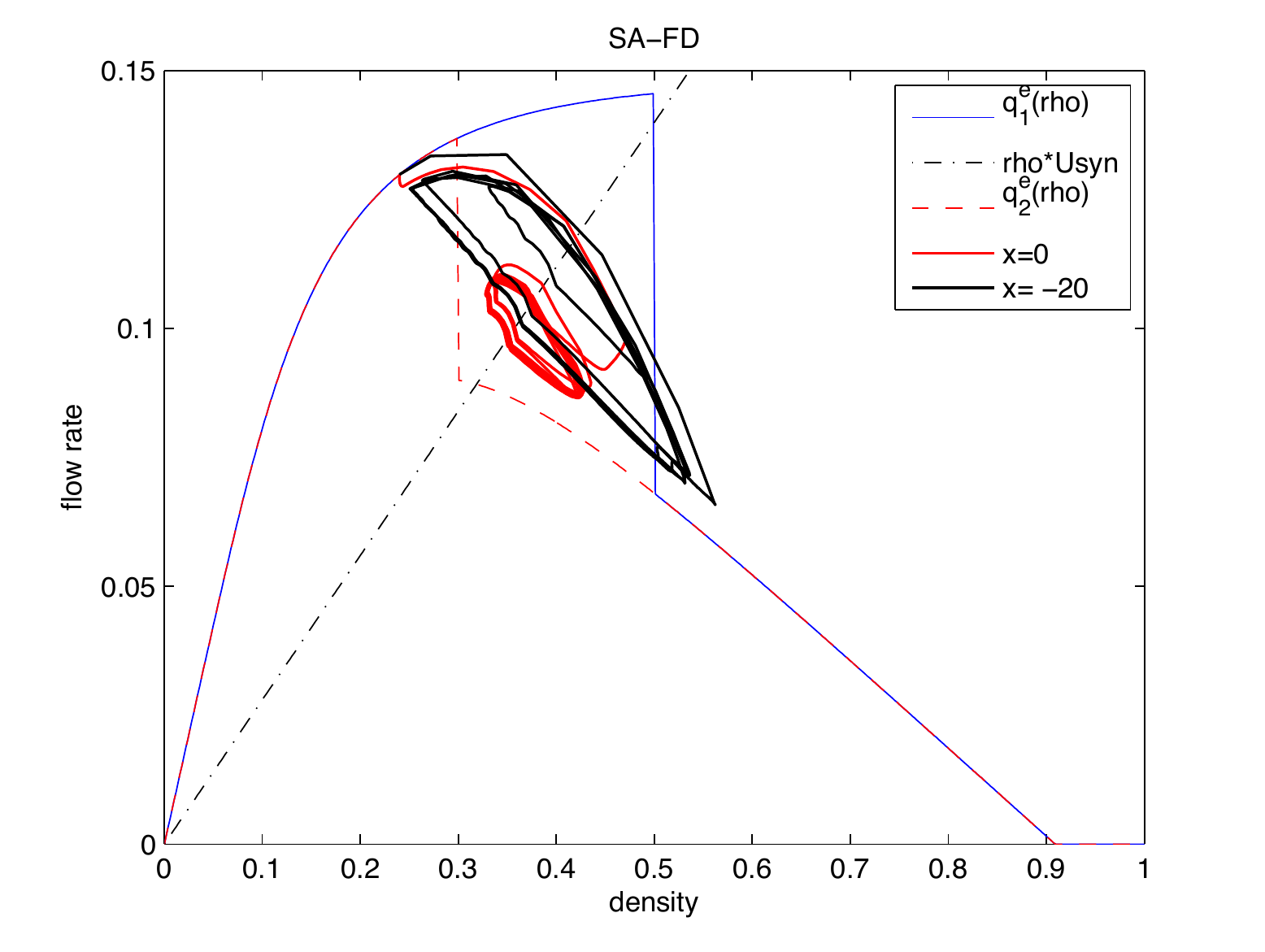,height=.25\textheight,width=.45\textwidth}
\epsfig{file=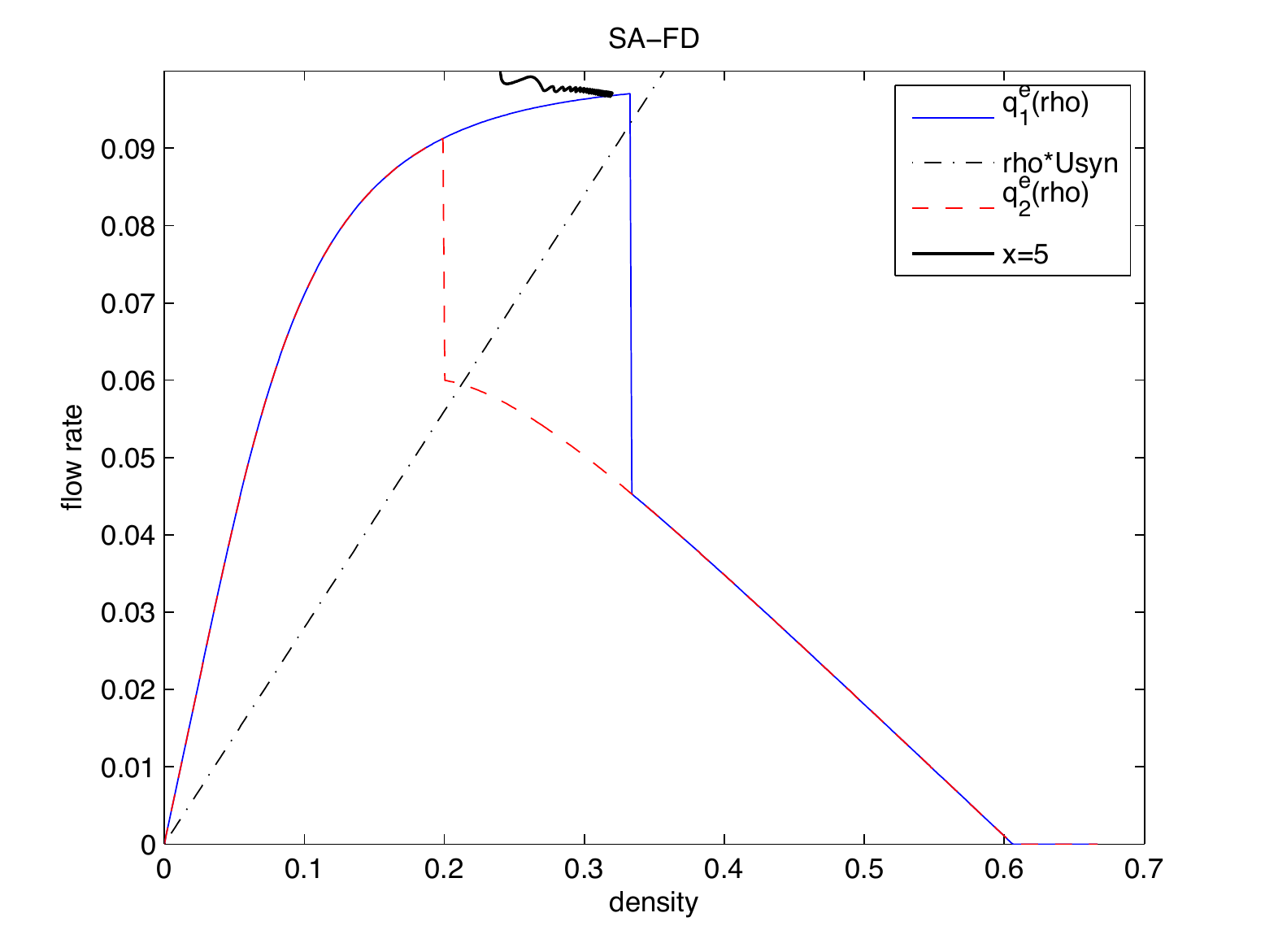,height=.25\textheight,width=.45\textwidth}
\epsfig{file=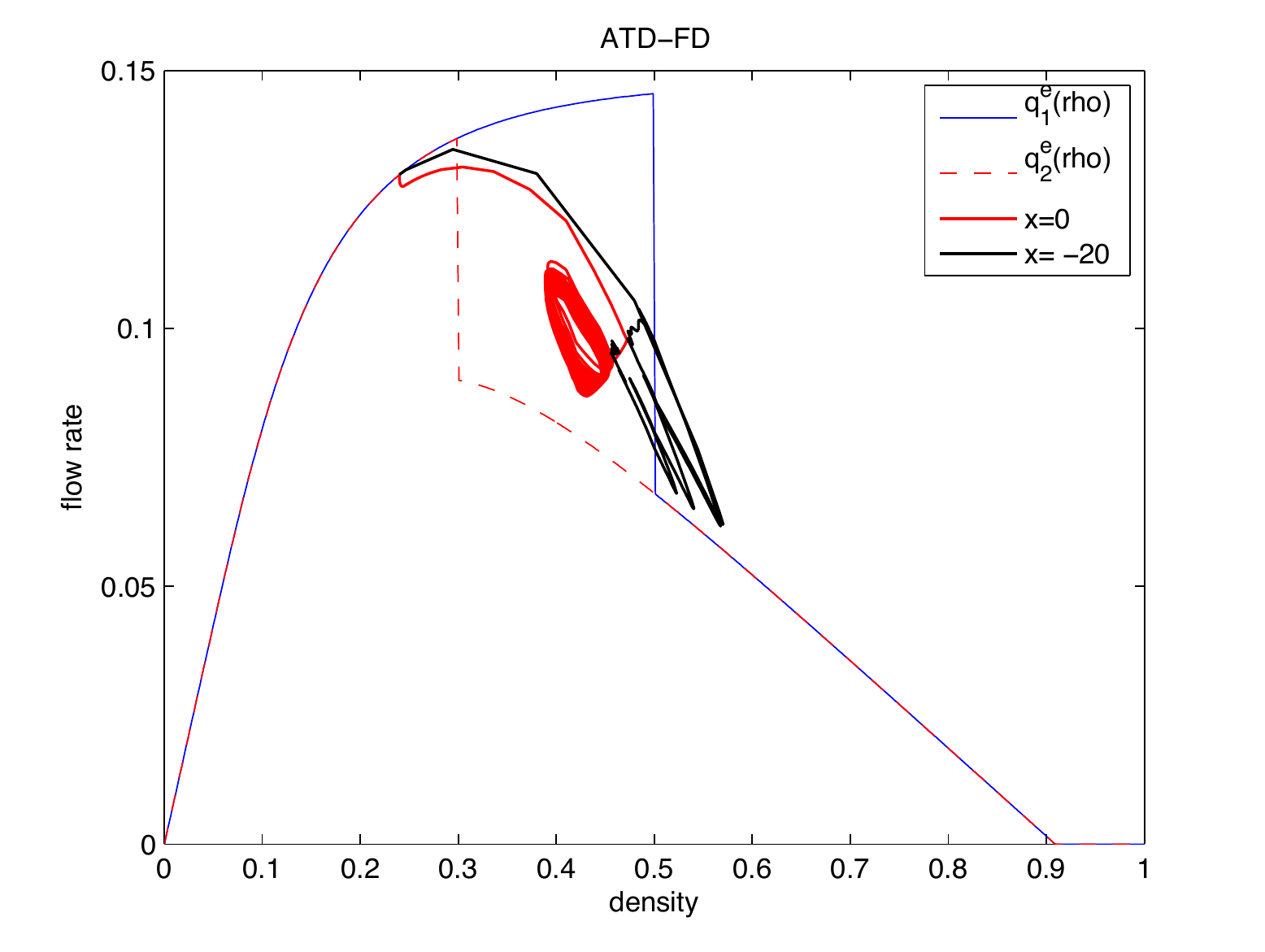,height=.25\textheight,width=.45\textwidth}
\epsfig{file=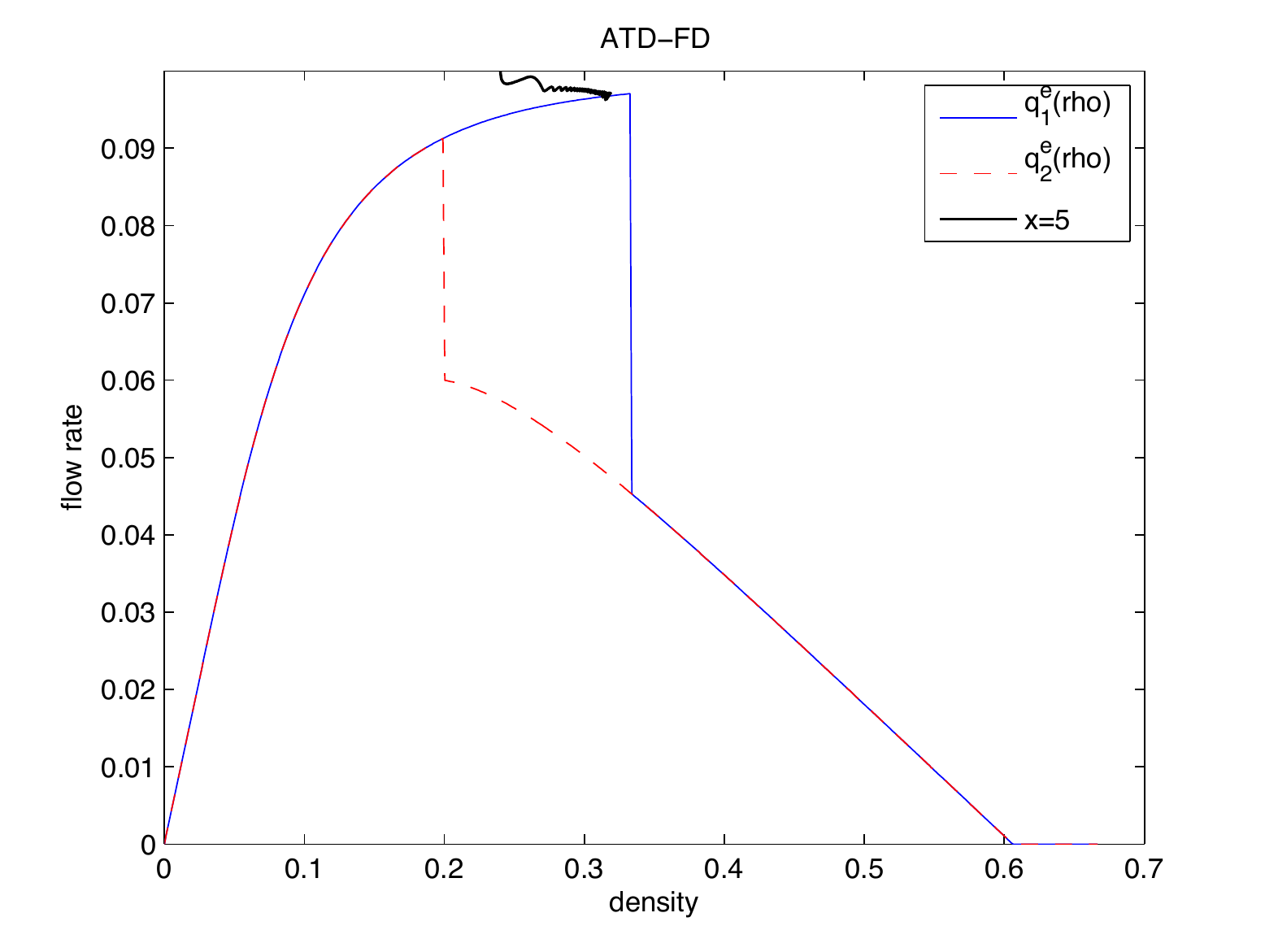,height=.25\textheight,width=.45\textwidth}
\end{center}
\caption{$ \rho u-\rho$ relation for the kinetic, switching curve, SA and ATD models upstream ($x=-20$), within ($x=0$) and downstream ($x=5$) of the bottleneck.}
\label{fdrelation}
\end{figure}

Finally, the macroscopic equations are investigated for a bottleneck
situation. For the computations we choose a Godunov method, see \cite{Tor}. We use a mesh size $\Delta x=0.15$, the Courant number $\lambda_{cfl}=0.99$ and a computation time $T_{end}=400$ units.
Figures \ref{rhofluxkinetic}, \ref{rhofluxswitch}, \ref{rhofluxSA} and \ref{rhofluxATD} show the velocity and density in space and time for a three lane highway with a reduction of lanes from $3$ to $2$ at $X=0$ for the four different models. In the simulation, the lane reduction is achieved by multiplying the density $\rho$ in the term on the right hand side of the equations by a factor $\frac{2}{3}$ for $X>0$ units and using a linear  interpolation between the 2 regions.
Apart from the ATD-type model, one clearly observes large changes in velocity and density in the backwards propagating traffic jams which might be interpreted as stop and go behaviour.  Figure \ref{fdrelation} show the flow-density relation at various locations of the considered highway, i.e. upstream of the bottleneck ($X=-20$), within the bottleneck ($X=0$),and downstream of the bottleneck ($X=5$). The flow rate drops from the initial value used in the simulation to settle at the maximum values for the considered highway's downstream location, $X=5$.  
For the ATD model in its present form we obtain a rather different behavior  due to the zero forcing inside the multi-valued region.

\begin{remark}
The models differ in the frequency and uniformity of the waves generated at the bottleneck. 
The fact that the ATD model does not generate stable waves in this situation does not mean that the model is in general incapable to 
describe traffic situation with such patterns. For example, in \cite{Ker03} situations are described where these waves appear. 
The models derived from the kinetic equations can be viewed as intermediate models between switching curve and ATD model, compare 
Figure \ref{uekineticcut}.
\end{remark}

\begin{remark}
A similar investigation could be performed for the model in \cite{BWGPB,Col} with suitable right hand side.
\end{remark}

\begin{remark}
Using the coefficient $c(\rho) =C \rho$  as done for example in \cite{CG07}, we obtain similar simulation results as above, if the parameters
 are suitably choosen. 
\end{remark}

\begin{remark}
The stable waves excited by small periodic perturbations as discussed in \cite{GKR03, Gre04} which may also be obtained from equations with single valued right hand sides are usually not persistent anymore for bottleneck situations. These waves are damped out as the high density region travels backward from the bottleneck.
\end{remark}

\subsubsection*{Summary}

Multi-valued fundamental diagrams are obtained using different approaches: a derivation from microscopic equations given in
\cite{Ker99,Ker03}, from  kinetic models as in \cite{GKMW03} and 
a phenomenological macroscopic model from \cite{GKR03}. 
These approaches are 
compared with each other from the point of view of their multi-valued fundamental diagrams and for an inhomogeneous bottleneck simulations
without any external excitation.
Apart from the ATD-type model, all the other models are able to show stop and go patterns for the described situation
with a bottleneck without external excitation of waves by the ingoing flow from an  on-ramp.

\subsubsection*{Acknowledgments}

The present work has been supported by the DAAD, PhD-Program MIC, Kaiserslautern.

\bibliographystyle{siam}

\end{document}